\documentclass{article}
\usepackage{latexsym}
\usepackage{amsfonts}
\usepackage{epsfig}
\usepackage{stmaryrd}
\usepackage{pxfonts}
\usepackage{wasysym}
\usepackage{bm}
\usepackage{algorithm}
\usepackage{algorithmic}
\usepackage[english]{babel}
\usepackage{subfigure, booktabs}

\title{A comparison of higher-order weak numerical schemes for stopped stochastic differential equations}
\date{\vspace{-5ex}}
 \author{ Francisco Bernal \footnotemark[1] \ \footnotemark[3]
\and Juan A. Acebr\'{o}n \footnotemark[2]\ \footnotemark[1]}

\begin{document}
\maketitle

\renewcommand{\thefootnote}{\fnsymbol{footnote}}

\footnotetext[1]{INESC-ID$\backslash$IST, TU Lisbon. Rua Alves Redol 9, 1000-029 Lisbon, Portugal.}
\footnotetext[2]{ISCTE - Instituto Universit\'ario de Lisboa
Departamento de Ci\^{e}ncias e Tecnologias de Informa\c{c}\~{a}o. Av. das For\c{c}as Armadas 1649-026 Lisbon, Portugal. ({\tt juan.acebron@ist.utl.pt})}
\footnotetext[3]{Center for Mathematics and its Applications,
Department of Mathematics, Instituto Superior T\'ecnico. 
Av. Rovisco Pais 1049-001 Lisbon, Portugal. ({\tt francisco.bernal@ist.utl.pt})}
\renewcommand{\thefootnote}{\arabic{footnote}}

\begin{abstract}
We review, implement, and compare numerical integration schemes for spatially bounded diffusions stopped at the boundary which possess a convergence rate of the discretization error with respect to the timestep $h$ higher than ${\cal O}(\sqrt{h})$. We address specific implementation issues of the most general-purpose of such schemes. They have been coded into a single Matlab program and compared, according to their accuracy and computational cost, on a wide range of problems in up to ${\mathbb R}^{48}$. The paper is self-contained and the code will be made freely downloadable.
\end{abstract}

{\em Keywords:} Weak convergence, Feynman-Kac, stochastic differential equation, bounded diffusion, first-exit problem.

\section{Introduction}\label{S:Introduction}
Producing numerical approximations, via Monte Carlo simulations, to the expected value of functionals involving stochastic differential equations (SDEs) stopped at a boundary is pervasive in scientific computing. Usually, the quantities of interest are the mean first passage time of a walker from that domain \cite{Redner2001} or the pointwise solution of a boundary value problem (BVP) \cite{Freidlin_Book}. Some relevant fields where it encounters application are: biology \cite{Tamborrino2014}, barrier problems in finance \cite{Glasserman2003}, mixing of fluids \cite{Iyer2010} and analysis of noisy dynamical systems \cite{Schwabedal_PRL}--where the interest often lies on first exit times--, as well as medical imaging \cite{Maire&Simon}, chemistry \cite{Mascagni2002}, design of integrated circuits \cite{Bernal2014}, and high-performance supercomputing \cite{Acebron2005}--more focused on the connection to BVPs. It is well known that the convergence rate with respect to the timestep $h$ of the approximation to the expected value--the weak convergence rate of the numerical scheme--is only ${\cal O}(\sqrt{h})$, unless specific measures to handle the interaction between the diffusion random path and the boundary are undertaken \cite{Gobet2001_BB}. 
This compares very unfavourably with the case of diffusions in free space, where even the most basic method--namely Euler-Maruyama's--has a weak order ${\cal O}(h)$ \cite{Kloeden&Platten}. In the last two decades, a number of weak schemes for bounded diffusions have been put forward seeking to raise the weak convergence rate to ${\cal O}(h^\delta)$, with $1/2<\delta\leq 1$. Henceforth, we will refer to a scheme with $\delta> 1/2$ as a "higher-order scheme". (We will use both terms "scheme" and "integrator" interchangeably.)

In applications, such a low exponent as $1/2$ means that in order to bound the total error by a tolerance $a$, a very small $h$ is needed, which leads to a cost ${\cal O}(a^{-4})$ \cite{Giles_seminal} and so to unpractically lengthy simulations if good accuracy is sought for. This exponent stems from the need to balance the statistical (Monte Carlo) error with the weak order of convergence of the integrator (i.e. the bias). In fact, the exponent can be split as $2+1/\delta$, where $2$ and $1/\delta$ are the contribution of either factor. Consequently, raising $\delta$ to $1$ means that the cost drops to ${\cal O}(a^{-3})$. In order to reduce it further, a suitable variance reduction method (tackling the $'2'$ in the exponent) and/or a $\delta>1$ are required. 

Giles' Multilevel method \cite{Giles_seminal} is the algorithm of choice for variance reduction, and has truly revolutionized the field of stochastic simulation. In it, several versions of the same stochastic paths--solved at ever smaller values of $h$--are combined into a telescopic estimator of the expected value in such a way that the rougher versions serve as control variates for the finer ones, thus reducing the overall variance. Importantly, within the Multilevel loop the trajectories are numerically integrated with a legacy scheme such as those considered in this article. The Multilevel method was first extended to stopped diffusions in \cite{Higham2013} using the Euler-Maruyama scheme as integrator, achieving a simulation cost of ${\cal O}\big(a^{-3}|\log{a}|^{1/2}\big)$. This has been recently reduced down to ${\cal O}\big(a^{-2}|\log{a}|^3\big)$ in \cite{Giles_y_yo}, by improving the convergence rate of the variance to linear and replacing the Euler-Maruyama  integrator by the Gobet-Menozzi one. The difficulty of combining Multilevel with a given integrator lies in the analysis of the variance, which depends on the {\em strong} convergence properties of the integrator. In any case, the {\em weak} order of convergence of the Multilevel estimator is that of the integrator inside it. Moreover, the Multilevel error estimate (see \cite{Giles_seminal}) also relies on that weak order being known in advance--as it is the case whenever the SDE is to be solved to within a given tolerance $a$.  

For unbounded diffusions, there exist superlinear weak integrators \cite{Kloeden&Platten}. Moreover, Talay-Tubaro extrapolation can be used \cite{Talay&Tubaro90}, also combined with Multilevel \cite{Giles_seminal,Lemaire&Pages}. On the other hand, while there are no superlinear integrators for bounded SDEs (to the authors' best knowledge),  extrapolation could in principle still be used \cite[chapter 5]{Milstein_Tretyakov_Book}. As a related idea, linear regression was used in \cite{LTMCR} to substantially reduce the bias of bounded SDEs. But again, all those strategies critically depend on the prior knowledge of the integrator's order $\delta$. 

Given the importance of this topic, it may be surprising that no systematic comparison of such higher-order schemes has been carried out before (as far as we know). This is the main purpose of this paper. A second aim is to 
identify, and hopefully clarify, the specific implementation issues with each specific integrator.
   
For the sake of conciseness, this paper is restricted to autonomous bounded diffusions stopped at the boundary.
(We recall that, with killed diffusions, one is interested in finding out whether or not the process has hit the boundary--typically, a barrier--within some time interval. On the other hand, with stopped diffusions, the exact hitting point is relevant, too.) In Section (\ref{S:Representation}), this problem will be explicitly connected with the pointwise solution of elliptic boundary value problems (BVPs) with Dirichlet boundary conditions (BCs). In general, the schemes that we review here can very easily be adapted to time-dependent, bounded SDEs (related, in turn, to parabolic BVPs with Dirichlet BCs). On the other hand, bounded SDEs reflected on the boundary require the approximation of the local time on the reflecting boundary, for which a different machinery is required \cite{Freidlin_Book,Constantini1998,Gobet2001_BB}. (For the sake of completeness, we just point out that those SDEs are related to BVPs with Neumann and Robin BCs.) 

Moreover, we have chosen to review only higher-order integrators which are of {\em general purpose}. By this, we mean that they are not restricted, by construction, to particular regular geometries or SDEs. For instance, the method of choice to solve Laplace's equation with Monte Carlo is the Walk on Spheres \cite{Muller1956}, which takes full advantage of the isotropy of the diffusion. Particularizing to convex polyhedric domains, the Walk on Rectangles algorithm \cite{Deaconu2006} is preferable, for it deals efficiently with corners. These two methods can be extended to slightly more complex SDEs, such as those having gradient and potential terms with constant coefficients. Other methods for isotropic SDEs are: exponential timestepping \cite{Jansons_Lythe_1st}, designed for parabolic boundaries; and the adaptive method in \cite{Bayer2010}, which may be the best option if they have localized features such as boundary singularities. Finally, Mair\'e and his collaborators have developed an accurate technique for dealing with the boundaries based on "randomization" of the BCs \cite{Maire2008}, which has been successfully used in applications, like \cite{Maire&Simon}. Randomization can also be applied to transmission BCs \cite{Mascagni2002,Bernal2014}. Most of the previous methods can be adapted to reflecting BCs.

Based on the general-purpose criterion, we have included the following integrators into the survey: Brownian bridge \cite{Mannella1999}, Buchmann and Petersen's \cite{Buchmann_Petersen_SIAM,Buchmann_JCP}, Gobet and Menozzi's \cite{Gobet&Menozzi_2010}, and two methods by Milstein \cite{Milstein97,Milstein&Tretyakov_simplest}. The first three are based on the Euler-Maruyama method, which is also included as a reference. 
As a matter of fact, the method of Buchmann and Petersen is restricted to diagonal diffusions, but given its excellent performance with high-dimensional Poisson's equations, we have included it nonetheless.    

\begin{figure}[h]
\centerline{\includegraphics[width=1.\columnwidth]{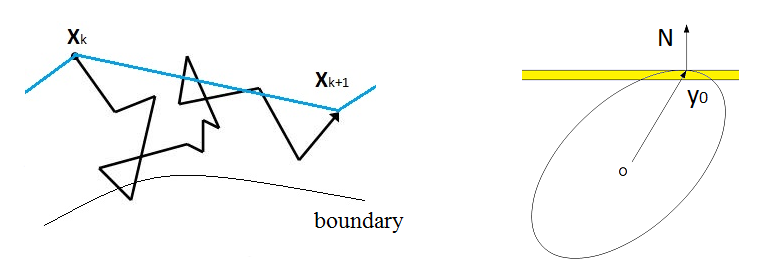}}
\caption{On the left, an excursion outside $\Omega$ missed by the coarser time discretization between ${\bf X}_k$ and ${\bf X}_{k+1}$. Right: tangent ellipsoid with a given orientation. The absorption layer (of width $\sqrt{Dh}$) is highlighted (see Section \ref{SS:WoE}).}
\label{I:Figura1}
\end{figure}

For stopped (as well as for killed) bounded diffusions, the reason for the lower weak rate of convergence compared to diffusions in free space is the ambiguity in determining the {\em first} intersection with a boundary of a random polygonal path--the numerical approximation to the diffusion under timestep $h$. See Figure \ref{I:Figura1} (left). The most natural way of approximating that first intersection is taking some point between the last iteration inside and the first interation outside. Unfortunately, not only does this result in a weak rate ${\cal O}(h^{1/2})$, but also leads to a systematic {\em overestimation} of the mean first exit time. In order to fix both issues, missed excursions such as in Figure \ref{I:Figura1} (left) must be taken into account. In the literature, three approaches have been proposed:
\begin{itemize}
\item Compute the probability that an excursion took place during every discrete step inside and terminate accordingly (used in the Brownian bridge and in the related method by Buchmann and Petersen). 
\item Make up for the missed first exits by shrinking the domain (boundary shift method by Gobet and Menozzi).
\item Control the step length at every step in order to make sure that the trajectory either stays within the domain, or touches its boundary exactly at one known point (methods of bounded increments by Milstein). 
\end{itemize}

In order to actually enjoy a higher order of convergence, the integrators need that the SDE coefficients or the integration domain meet some regularity conditions, which applications often lack. On the other hand, it is interesting to check how the various higher-order methods fare precisely in that situation. Therefore, we have made a point of including both kind of domains (those very smooth and those with corners) in our comparison.  

The remainder of the paper is organized as follows. In Section \ref{S:Representation}, the precise connection between autonomous SDEs stopped at the boundary and linear elliptic BVPs with Dirichlet BCs is made. We have put together the common aspects to all of the integrators in Section \ref{S:Preliminaries}. Then, the five methods singled out in the introduction are reviewed in Section \ref{S:Integrators}. The test problems are presented in Section \ref{S:Experiments}, and Section \ref{S:Conclusions} concludes the paper. 

\section{Probabilistic representation}\label{S:Representation}
In this paper we consider expectations of the form

\begin{equation}
\label{F:Dynkin}
u({\bf x}_0)= {\mathbb E}\Big[\, g({\bf X}_{\tau})e^{\int_0^{\tau}c({\bf X}_s)ds} \, + \, \int_0^{\tau}f({\bf X}_t)e^{\int_0^t c({\bf X}_s)ds}dt\,\big|\,{\bf X}_0={\bf x}_0\,\Big],
\end{equation}

where ${\bf X}_{\tau}$ is the value at $t=\tau>0$ of the stochastic process ${\bf X}_t:[0,\tau]\rightarrow\mathbb{R}^D$, driven by the SDE 

\begin{equation}
\label{F:SDE}
d{\bf X}_t= {\bf b}({\bf X}_t)dt + \sigma({\bf X}_t)d{\bf W} _t\,\qquad{\bf X}_0={\bf x}_0.
\end{equation}

In (\ref{F:SDE}), $\Omega\subset{\mathbb R}^D,\, D\geq 1$ is a bounded domain, $\sigma$ is the diffusion matrix, and ${\bf W}_t$ is a standard $D$-dimensional Wiener process. The {\em first exit time} (or {\em first passage time}) $\tau$ is defined as 
\begin{equation}
\label{F:tau_def}
\tau= \inf_{t\geq 0} {\bf X}_t\in\partial\Omega
\end{equation}  
{\em i.e.}, the time when a solution of the SDE (\ref{F:SDE}) first touches $\partial\Omega$ at the {\em first exit point} ${\bf X}_{\tau}$. The process ${\bf X}_t,\,0\leq t\leq \tau$ can be thought of as a non-differentiable trajectory inside $\Omega$. 

If $c=g=0$ and $f=1$, (\ref{F:Dynkin}) is the mean exit time from $\Omega$ of a diffusion starting at ${\bf x}_0\in\Omega$. The other most important case takes place when $\det({\sigma})\neq 0$, $c<0$ and
the coefficients $a_{11},\ldots,a_{DD},b_1,\ldots,b_D,c,f,g:{\mathbb R}\mapsto{\mathbb R}$ are assumed regular enough that the solution to the elliptic BVP with Dirichlet BCs (\ref{F:EllipticBVP}) exists and is unique \cite{Miranda70}:
\begin{equation}
\label{F:EllipticBVP} 
\sum_{i=1}^D\sum_{j=1}^D a_{ij}\frac{\partial^2 u}{\partial x_i\partial x_j} + \sum_{k=1}^D b_k\frac{\partial u}{\partial x_k} + cu + f = 0, 
\textrm{ if } {\bf x}\in\Omega;\qquad
u({\bf x})= g, \textrm{ if } {\bf x}\in\partial\Omega,
\end{equation}

where the matrix 
$A({\bf x})=[a_{ij}]=\sigma\sigma^T/2$ is positive definite. Then, (\ref{F:Dynkin}) is the pointwise solution of (\ref{F:EllipticBVP}) at ${\bf x}_0$. Formula (\ref{F:Dynkin}) is Dynkin's formula, a particular case of the more general Feynman-Kac formula for parabolic BVPs--see \cite{Freidlin_Book} for more details. In the remainder of the paper, we will assume that the identification between the SDE (\ref{F:SDE}) and the BVP (\ref{F:EllipticBVP}) holds and refer to (\ref{F:Dynkin}) both as the solution of an SDE or a BVP.

An alternative representation of $u({\bf x}_0)$, due to Milstein \cite{Milstein_Tretyakov_Book}, is

\begin{equation}
\label{F:gXYplusZ}
u({\bf x}_0)= {\mathbb E}[\,\phi({\bf x}_0)\,]:= {\mathbb E}\Big[\, g({\bf X}_{\tau})Y_{\tau} + Z_{\tau} \,\Big], 
\end{equation}

where $\phi({\bf x})$ are called the {\em scores} and $({\bf X}_{\tau},Y_{\tau},Z_{\tau})$ are the evaluation at time $t=\tau$ of the solution of the system of SDEs 

\begin{equation}
\label{F:Milstein_sys}
\left\{
\begin{array}{lr}
d{\bf X}_t= {\bf b}({\bf X}_t)dt + \sigma({\bf X}_t)d{\bf W}_t, & {\bf X}_0={\bf x}_0, \\
dY= c({\bf X}_t)Ydt, & Y(0)=1,\\
dZ= f({\bf X}_t)Ydt, & Z(0)= 0.
\end{array}
\right. 
\end{equation}

Numerically, the expected value (\ref{F:gXYplusZ}) is approximated with a Monte Carlo (MC) method involving $N>>1$ independent realizations of (\ref{F:Milstein_sys}), which are in turn integrated numerically with a suitable scheme (integrator) based on a small timestep $h>0$, which we will consider constant:
\begin{equation}
u({\bf x}_0)\approx u_{h,N}({\bf x}_0):=\frac{1}{N}\sum_{j=1}^N \phi_h^{(j)}({\bf x}_0),
\end{equation}
where $\phi_h^{(j)}({\bf x}_0)$ is the $j^{th}$ realization of (\ref{F:Milstein_sys}) with a timestep $h$. In doing this, two kinds of numerical errors arise: the time-stepping error involved in the numerical integration (called bias), and the statistical error, due to replacing the mean with the expected value. 
The bias depends on the integrator. By virtue of the central limit theorem, the statistical error is $\sqrt{V[\phi({\bf x}_0)]/N}$, where $V[\phi({\bf x}_0)]$ is the variance of the scores \cite{Kloeden&Platten}. Extending well-established results for the case of unbounded diffusions \cite{Giles_seminal}, the root mean square (RMS) error of the MC approximation to (\ref{F:EllipticBVP}) as $h\shortrightarrow 0^+,N\shortrightarrow \infty$ obeys
\begin{equation}
\label{F:MC_errors}
\epsilonup:=\sqrt{RMS\big(|u({\bf x}_0)-u_{h,N}({\bf x}_0)|\big)}\leq  |C|h^{\delta} + q\sqrt{\frac{V[\phi({\bf x}_0)]}{N}}\textrm{ with probability $P_q$,}
\end{equation} 
where $C$ is a constant, $\delta>0$; and $P_q\approx 68.3\%, 95.5\%$, and $99.7\%,$ for $q=1,2,3$. 
 The value of $\delta$ for different integrators is the object of study of this paper.

\subsection{Alternative representation}\label{SS:VR}

Consider the related system to (\ref{F:Milstein_sys}): 

\begin{equation}
\label{F:VR_Milstein_sys}
\left\{
\begin{array}{lr}
d{\tilde {\bf X}}_t= [b({\tilde {\bf X}}_t)-\sigma({\tilde {\bf X}}_t){\bm \mu}({\tilde {\bf X}}_t)]dt + \sigma({\tilde {\bf X}}_t)d{\bf W}_t, & {\tilde {\bf X}}_0={\bf x}_0, \\
d{\tilde Y}= c({\tilde {\bf X}}_t){\tilde Y}dt + {\bm \mu}^T({\tilde {\bf X}}_t){\tilde Y}d{\bf W}_t, & {\tilde Y}(0)=1,\\
d{\tilde Z}= f({\tilde {\bf X}}_t){\tilde Y}dt + {\bf F}^T({\tilde {\bf X}}_t){\tilde Y}d{\bf W}_t, & {\tilde Z}(0)= 0,
\end{array}
\right. 
\end{equation}

and let ${\tilde \phi}({\bf x}_0)=g({\bf {\tilde X}}_{\tau}){\tilde Y}_{\tau}+{\tilde Z}_{\tau}$. It can be proved \cite[chapter 2]{Milstein_Tretyakov_Book} that
${\mathbb E}[{\tilde\phi}({\bf x}_0)]={\mathbb E}[\mathbb\phi]=u({\bf x}_0)$, but $V[{\tilde\phi}({\bf x}_0)]\neq V[\phi({\bf x}_0)]$. 

We will exploit (\ref{F:VR_Milstein_sys}) in two ways. First, it allows us to remove the drift, which is crucial in the Walk on Ellipsoids method in Section \ref{SS:WoE}. Secondly, if ${\bm\mu}({\bf x})$ and ${\bf F}({\bf x})$ are so chosen that
\begin{equation}
\label{F:Pathwise_VR}
\sigma^T\nabla u+u{\bm\mu}+{\bf F} =0,
\end{equation}
then $V[\tilde\phi({\bf x}_0)]=0$. Unfortunately, this requires the solution of (\ref{F:EllipticBVP}) in the first place. However, it will allow us to accelerate the simulations in Section \ref{S:Experiments} (because far fewer trajectories are required for a given statistical error), and so to explore deeper the weak error curves of the various integrators.

\section{Preliminaries}\label{S:Preliminaries}

\subsection{Boundary data} In this paper, $||\cdot||$ is always the 2-norm. The signed distance from a point ${\bf x}$ (in- or outside $\Omega$) to the boundary is denoted as $d({\bf x})$, with the criterion that $d({\bf x})$ is negative inside $\Omega$, zero on $\partial\Omega$, and positive outside. The projection of  ${\bf x}$ onto the boundary $\partial\Omega$ is the closest point to ${\bf x}$ on $\partial\Omega$ and is denoted by $\Pi_{\partial\Omega}({\bf x})$. This projection is considered unique, which may not be true unless $|d({\bf x})|=||{\bf x}-\Pi_{\partial\Omega}({\bf x})||$ is small enough. $\partial\Omega$ is assumed smooth enough that the outward unit normal is uniquely defined on any point of it, so that there is a well-defined tangent hyperplane. The notation ${\bf N}(\bf x)$ and $H_{\partial\Omega}({\bf x})$ stand for the normal and the tangent hyperplane evaluated at $\Pi_{\partial\Omega}({\bf x})$. 

\subsection{Signed distance map} For a given point ${\bf x}\in{\mathbb R}^D$, all the numerical integrators need to know the {\em boundary data} $d({\bf x}),\Pi_{\partial\Omega}({\bf x}),{\bf N}({\bf x})$, and $H_{\partial\Omega}({\bf x})$. Unless $\Omega$ has a very simple geometry (or is a combination thereof), the signed distance $d(\bf x)$ cannot be computed analytically, but the eikonal equation must be solved in advance \cite{Bernal2014}:
\begin{equation}
\label{F:Eikonal}
||\nabla d({\bf x})||=1\textrm{ if ${\bf x}\notin\partial\Omega$},\qquad d({\bf x})=0\textrm{ if ${\bf x}\in\partial\Omega$.}
\end{equation}
Once the signed distance map is available, ${\bf N}({\bf x})=\nabla d({\bf x})$ (hence the sign convention), $\Pi_{\partial\Omega}({\bf x})={\bf x}-d({\bf x})\nabla d({\bf x})$, and $H_{\partial\Omega}({\bf x})=\{{\bf y}\,|\,\big({\bf y}-\Pi_{\partial\Omega}({\bf x})\big)^T{\bf N}({\bf x})=0\}$. The numerical procedure of choice for solving (\ref{F:Eikonal}) is the Fast Marching method. For more examples and details on Matlab implementations in the context of SDEs see \cite{Bernal2014}. We remark that, since the boundary data are needed only to handle the interaction trajectory-boundary, the signed distance is usually required in a shell of ${\cal O}(\sqrt{h})$ thickness around $\partial\Omega$, not in the whole of ${\mathbb R}^D$. 

\subsection{Timestepping notation} Let $k=0,1,2,\ldots$ be the {\em steps} and $h>0$ the timestep. The numerical approximation to ${\bf X}_t,\,t=kh$ is ${\bf X}_k$ (and analogously for $Y_k$ and $Z_k$). The notation $\sigma_k$ stands for $\sigma({\bf X}_k)$; ${\bf b}_k,c_k$ and $f_k$, as well as $d_k$ and ${\bf N}_k$ are defined in the same way. The numerical approximation to the exit time $\tau$ is represented by $\nu$. Finally, $\phi_h:=g({\bf X}_{\nu})Y_{\nu}+Z_{\nu}$ are the numerical approximations to $\phi,{\bf X}_{\tau},Y_{\tau}$, and $Z_{\tau}$, respectively. 
 
\subsection{\bf Probability distributions and integrators} The notation $w\sim P$ means that the stochastic variable $w$ has been drawn from the probability distribution $P$. Table \ref{T:Distributions} displays the probability distributions quoted in the description of the integrators. The initials EM, GM, BB, BP, RW, and WoE stand for the Euler-Maruyama, Gobet-Menozzi, Brownian bridge, Buchmann-Petersen, Milstein's random walk and Walk on Ellipsoids integrators, respectively.

\begin{table}[H]
\begin{footnotesize}
\caption{Probability distributions used in this paper.}
\label{T:Distributions}
\[\begin{array}{llll}
\hline\noalign{\smallskip}
\textrm{Symbol} &  \textrm{Name} & \textrm{Integrator} & \textrm{Matlab command}\\
\noalign{\smallskip}\hline\noalign{\smallskip}
{\cal U}&\textrm{Bernoulli (uniform in $[0,1]$)} & BB,BP,RW &rand\\
{\cal N}& \textrm{standard normal} & EM,GM,BB,BP & randn\\
{\cal IG}(\gamma,\delta)& \textrm{inverse Gaussian} & BP &igpdf\textrm{ (see text)}\\
{\cal B} & \textrm{binary ($\pm 1$ with $50\%$ chance})& RW & 2*randi-1 \\
{\cal S}(D) & \textrm{uniform on the surface of a unit $D$-sphere}& WoE & dboluda\textrm{ (see text)}\\
\hline
\end{array}\]
\end{footnotesize}
\end{table}

The distribution ${\cal S}(D)$ can be constructed according to \cite[remark 5.1.3]{Milstein_Tretyakov_Book}, but we provide a Matlab code based on a simpler idea below, as well as for the inverse Gaussian distribution (see also \cite{Michael1976}).
\begin{footnotesize}
\begin{verbatim}
function w= dboluda(N,D) 
x= rand(N,D)-1/2; w= w./sqrt(sum(x.^2,2)); 
\end{verbatim}

\begin{verbatim}
function [res]= igpdf(gamma,delta,sampleSize)
chisq1 = randn(sampleSize).^2;
res = delta + 0.5*delta./gamma .* ( delta.*chisq1 - sqrt(4*delta.*gamma.*chisq1 + (delta.^2).*chisq1.^2) );
l = rand(sampleSize) >= delta./(delta+res);
res( l ) = (delta( l ).^2)./res( l );
\end{verbatim}
\end{footnotesize}

\subsection{SDE coefficients}\label{SS:Gershgorin} 
The functions $\sigma,{\bf b},c\leq 0,f$ and $g$ are assumed to have at least the regularity required for the Feynman-Kac formula to hold, and the matrix $A$ in (\ref{F:EllipticBVP}) is assumed symmetric positive definite in $\Omega$ (strict ellipticity condition). Thus, for any ${\bf v}\in{\mathbb R}^D$, $||\sigma{\bf v}||={\bf v}^T\sigma^T\sigma{\bf v}>0$. Then, $A$ has a Cholesky decomposition so that $\sigma$ can always be chosen lower triangular. Moreover, the eigenvalues of $A$ 
are $\lambda_{max}(A)\geq\lambda_2\geq\ldots\geq\lambda_D>0$ (note that they may depend on ${\bf x}$). Some integrators need $\lambda_{max}(A)$. Computing it may be very costly (especially in high dimensions, unless $A$ has an exploitable pattern), so that it may be worth bounding it above via Gershgorin's circle theorem: let $M$ be a square matrix and $\lambda$ any of its eigenvalues. Then,  
\begin{equation}
|\lambda-M_{ii}|\leq \sum_{i\neq j}|M_{ij}|. 
\end{equation}  
If $\lambda>0$ as with $A$, $|\lambda-M_{ii}|\geq \lambda-|M_{ii}|$ so that
\begin{equation}
\label{F:Gershgorin}
\lambda_{max}(A)\leq \max_{j=1,\ldots,D}\sum_{i=1}^D|a_{ij}|.
\end{equation}

\section{Review of the integrators}\label{S:Integrators}

\subsection{Boundary shift method of Gobet and Menozzi (GM)}\label{SS:GM}

The Euler-Maruyama method (EM) for (\ref{F:Milstein_sys}) is Algorithm \ref{A:EM}. As already pointed out, neglecting intermediate excursions leads to $0<{\mathbb E}\nu-{\mathbb E}\tau={\cal O}(\sqrt{h})$. Gobet and Menozzi's idea is to compensate for them by {\em shrinking} the domain: if one lets $\Omega_h'$ and $\partial\Omega_h'$ be the shrunken domain and its boundary, the {\em discrete} trajectories will {\em on average} hit $\partial\Omega_h'$ before they would hit $\partial\Omega$. Gobet and Menozzi prove that if the shift inwards is locally given by
\begin{equation}
\label{F:GM_shift}
\textrm{shift= $-0.5826||\sigma({\bf x}){\bf N}({\bf x})||\sqrt{h}{\bf N}({\bf x}),\qquad {\bf x}\in\partial\Omega$,}
\end{equation} 
then the weak convergence rate of the Euler-Maruyama scheme improves to {\em o}($\sqrt{h}$)--see \cite{Gobet&Menozzi_2010} for the technical conditions. 
Importantly, in many instances the weak convergence rate turns out to be in fact linear, although precise conditions for this were not yet determined in \cite{Gobet&Menozzi_2010}.

Based on this result, the Gobet-Menozzi scheme (GM) boils down to replacing line 3 in Algorithm \ref{A:EM} with "{\bf if }$d_k>-0.5826||\sigma_k{\bf N}_k||\sqrt{h}$", as well as line 1 with "$d_0<-0.5826||\sigma_0{\bf N}_0||\sqrt{h}$".    

\begin{algorithm}[h!]
\caption{Euler-Maruyama (EM)}
\begin{algorithmic}[1]
\STATE{{\bf Data:} $h>0$, ${\bf X}_0={\bf x}_0\in\Omega,{Y_0}=1,Z_0=0$, $(d_0<0,{\bf N}_0)$}
\FOR{$k=0,1,2,\ldots$ until EXIT}
	\IF{$d_k<0$}
		\STATE{Let ${\bm \omegaup}_{k+1} \sim {\cal N}$, and take one unbounded step according to:}
		\begin{eqnarray}
		\label{F:E_M}
		\left\{\begin{array}{l}
		{\bf X}_{k+1}= {\bf X}_k + h{\bf b}_k + \sqrt{h}\sigma_k{\bm \omegaup_{k+1}}, \\
		Y_{k+1}= Y_k + hY_kc_k,  \\
		Z_{k+1}= Z_k + hY_kf_k. \\
		\end{array}\right.
		\end{eqnarray}
		\STATE{Compute $(d_{k+1},{\bf N}_{k+1})$ according to the distance map of $\Omega$}
	\ELSE{}
		\STATE{EXIT: {\bf goto} 10}
	\ENDIF
\ENDFOR
\STATE{ Take $\nu=kh,{\bf X}_{\tau}=\Pi_{\partial\Omega}({\bf X}_k),\,Y_{\tau}=Y_k,\,Z_{\tau}=Z_k$, and $\phi_h= g({\bf X}_{\nu})Y_{\nu}+Z_{\nu}$.}
\end{algorithmic}
\label{A:EM}
\end{algorithm}

\subsection{Brownian bridge (BB)}\label{SS:BB}

The introduction of Brownian bridge techniques in the numerical solution of SDEs representing elliptic problems with Dirichlet BCs is inspired by the earlier application to barrier problems in finance. A simple barrier problem is as follows. Let $X(t)$ be a one-dimensional process driven by a constant drift $\hat{b}$ and a constant volatility ${\hat\sigma}$:

\begin{equation}
\label{F:OneD}
dX(t)= {\hat b}dt + {\hat\sigma}dW_t,\qquad X(0)=y.
\end{equation}

Assume that after $t=h>0$, $z:=X(h)$ is known. Then, the process (\ref{F:OneD}) in $0\leq t\leq h$ conditioned to $X(h)=z$ is called a Brownian   
bridge pinned at $(0,y)$ and $(h,z)$. Let $y<B>z$ be a barrier. One is interested in the probability of the Brownian bridge hitting the barrier, {\em i.e.} taking the value $X(\tau')=B$ for some $0<\tau'<h$. This probability turns out to be independent of ${\hat b}$ and is \cite[section 3.1]{Buchmann_JCP}
\begin{equation}
\label{F:Tunneling}
P(\tau'<h)= \exp{\Big( -\frac{2(B-y)(B-z)}{h{\hat\sigma}^2}\Big)}.
\end{equation}

This suggests the incorporation into the Euler-Maruyama method of a {\em boundary test} for detecting excursions during a given time step \cite{Baldi1995}. Assume a one-dimensional autonomous SDE with non-constant coefficients $b(X)$ and $\sigma(X)$ in $\Omega:=(-\infty,x_B]$, and for some $k\geq 1$ let $X_k\in\Omega$ and $X_{k+1}\in\Omega$ ({\em i.e.} $X_k<x_B>X_{k+1}$). Further, assume that the trajectory during the time interval $kh<t<(k+1)h$ can be approximated by the Brownian bridge pinned at $X_k$ and $X_{k+1}$ with frozen diffusion $\sigma(X_k)$. (Such an approximation may not always be consistent: see \cite{Giraudo1999}.) 
Now let $w$ be a variable uniformly distributed in $[0,1]$. If $w<P$, where
\begin{equation}
\label{F:Mannella_P}
P= \exp{\Big( -\frac{2(x_B-X_{k})(x_B-X_{k+1})}{h\sigma^2({X_k})} \Big)},
\end{equation}
then the trajectory is deemed to have undergone an excursion, so that it is finished at the iteration $k+1$ with $\nu=(k+1)h$. A finite interval $\Omega:=[x_A,x_B]$ can also be tackled in this way provided that $|x_A-x_B|>>\sqrt{h}$, so that (\ref{F:Mannella_P}) approximately holds ({\em i.e.} if the probability of hitting the farther end can be neglected). These ideas were implemented by Mannella in \cite{Mannella1999}, who also improved on the zero-order frozen-coefficient approximation taken here. 

Higher dimensions can be tackled via the {\em half-space approximation}, which reduces the interaction trajectory-boundary to the one-dimensional case. Let us consider a diffusion with constant coefficients like (\ref{F:OneD}) taking place in ${\mathbb R}^D\,(D\geq 2)$; a hyperplane $H$ defining two half-spaces; and two points ${\bf x}_1$ and ${\bf x}_2$ in the same half-space. The probability that a Brownian bridge pinned on ${\bf X}(0)={\bf x}_1$ and ${\bf X}(h)={\bf x}_2$ crosses $H$ within $h$ is still given by (\ref{F:Tunneling}), provided that $(B-x)(B-y)$ is replaced by $d_H({\bf x}_1)d_H({\bf x}_2)$, where $d_H({\bf x})\geq 0$ is the distance between $\bf x$ and $H$. Therefore, the Brownian bridge test can be extended to any dimension if the boundary can be locally approximated by a hyperplane. Intuitively this occurs when $\partial\Omega$ is locally smooth and its radius of curvature is much larger than $\sqrt{h}$ and $d_H$. Importantly, this does not hold close to a corner or a cusp regardless of $h$. 

The Brownian bridge with half-space approximation for the Euler-Maruyama integration of (\ref{F:SDE}) then consists in substituting $\partial\Omega$ with $H_{\partial\Omega}({\bf X}_k)$ and considering the frozen-coefficient Brownian bridge
\begin{equation}
d{\bf X}= {\bf b}_kdt + \sigma_kd{\bf W}_t,\textrm{with ${\bf X}(kh)={\bf X}_k$ and ${\bf X}((k+1)h)={\bf X}_{k+1}$}
\end{equation}
after computing the Euler-Maruyama step ${\bf X}_{k+1}$. The probability of the particle exitting at the iteration $k+1$ is then \cite{Gobet2001_BB}

\begin{equation}
\label{F:Brownian_bridge}
p= \left\{ 
\begin{array}{ll}
1, & \textrm{ if ${\bf X}_{k+1}\notin\Omega$,}\\
\exp{\Big( -\frac{2d_kd_{k+1}}{h{\bf N}_k^T\sigma_k\sigma_k^T{\bf N}_k}\Big)} & \textrm{ otherwise.}
\end{array}
\right.
\end{equation}

Note that the distances in (\ref{F:Brownian_bridge}) are now signed, but $d_kd_{k+1}>0$. The pseudocode is given as Algorithm \ref{A:BB}. Assuming great smoothness ($C^5$) of the domain, the drift, and $\sigma$, Gobet proves that the half-space approximation has a linear weak order of convergence \cite{Gobet2001_BB}. (So that taking the exit at the iteration $k$ instead of $k+1$--or averaging between them--has an effect of order $h$ and does not alter the weak convergence rate.) However, note that the proof pertains to {\em killed}, rather than {\em stopped} diffusions. Therefore, it is valid in the estimation of mean exit times or barrier options. On the other hand, the Feynman-Kac functional (\ref{F:Dynkin}) (or equivalently (\ref{F:gXYplusZ})), depends on the point where the BC is evaluated (unless $g\equiv 0$, or if it is one end of a one-dimensional interval, {\em i.e.} $D=1$).

\begin{algorithm}[h!]
\caption{Brownian bridge (BB)}
\begin{algorithmic}[1]
\STATE{{\bf Data:} $h>0$, ${\bf X}_0={\bf x}_0\in\Omega,{Y_0}=1,Z_0=0$, and $(d_0<0,{\bf N}_0)$}
\FOR{$k=0,1,2,\ldots$ until EXIT}
\STATE{Get ${\bf X}_{k+1},Y_{k+1}$ and $Z_{k+1}$ according to the Euler-Maruyama scheme (\ref{F:E_M})}
\STATE{Compute $(d_{k+1},{\bf N}_{k+1})$ according to the distance map of $\Omega$}
\IF{$d_{k+1}\geq 0$} 
	\STATE{EXIT: {\bf goto} 14} 
\ELSE{} 
	\STATE{Let $w'\sim {\cal U}$}
	\IF{$\displaystyle{w'<\exp{\Big(-\frac{2d_kd_{k+1}}{{h\bf N}_k^T\sigma_k\sigma_k^T{\bf N}_k}\Big)}}$ (boundary test)} 
	\STATE{EXIT: {\bf goto} 14}
	\ENDIF
\ENDIF
\ENDFOR
\STATE{ Take $\nu=kh,{\bf X}_{\tau}=\Pi_{\partial\Omega}({\bf X}_k),\,Y_{\tau}=Y_k,\,Z_{\tau}=Z_k$, and $\phi_h= g({\bf X}_{\nu})Y_{\nu}+Z_{\nu}$.}
\end{algorithmic}
\label{A:BB}
\end{algorithm}

\subsection{The method of Buchmann and Petersen (BP)}\label{SS:BP}

This method is an improvement for the case of a class of stopped diffusions of the previous Brownian bridge test. The idea is to determine--in a statistical sense--{\em when} within the last time step of length $h$ the trajectory exitted the domain, and {\em where} on the tangent hyperplane the crossing took place. Consider again the one-dimensional diffusion (\ref{F:OneD}) and the same notation as in (\ref{F:Tunneling}). There are two cases: 1) $y<B\leq z$, when the trajectory has exitted for sure, and 2) $y<B>z$, when there may have been an excursion. Then, the associated exit times $\tau_1$ and $\tau_2$ are sampled from the distributions $P_1$ and $P_2$ below \cite{Buchmann_JCP}:
\begin{eqnarray}
P_1(\tau_1)= \bm{1}_{[\tau_1<h]}\frac{B-y}{\sqrt{2\pi {\hat\sigma}^2 \tau_1^3}}\sqrt{\frac{h}{h-\tau_1}}\exp{\Big[-\frac{1}{2{\hat \sigma}^2}\big( \frac{(z-B)^2}{h-\tau_1} -\frac{(z-y)^2}{h} + \frac{(B-y)^2}{\tau_1} \big)\Big]},
\end{eqnarray} 
\begin{eqnarray}
P_2(\tau_2)= \exp{\Big(-\frac{2(B-z)(B-y)}{{\hat\sigma}^2\tau_2}\Big)}.
\end{eqnarray}

In order to sample from those distributions, Buchmann shows that
\begin{equation}
\label{F:Tau1}
w\sim {\cal IG}\big(\, \frac{(B-y)^2}{h{\hat\sigma}^2},\frac{B-y}{z-B} \,\big)
\Rightarrow \tau_2= \frac{hw}{1+w}.
\end{equation}

(look up Table \ref{T:Distributions} for reference), and

\begin{eqnarray}
\label{F:Tau2}
w'\sim U \Rightarrow \tau_1= -\frac{2(B-y)(B-z)}{{\hat\sigma}^2\log{w'}}.
\end{eqnarray}

Contrary to the plain Brownian bridge method, now a more accurate estimate $kh<\nu<(k+1)h$ of $\tau$ is available, which can be used to improve on the exit point--where the Dirichlet BC is evaluated. In order to do so, Buchmann and Petersen use the half-space approximation to reason as follows \cite{Buchmann_Petersen_SIAM}. Assume that the trajectory is stopped at a time $kh+\tau'$ ($0<\tau'<h$) between ${\bf y}={\bf X}_k$ and ${\bf z}={\bf X}_{k+1}$ (either because ${\bf z}\notin\Omega$ or because an excursion was detected) at some exit point ${\bf y}+{\bm \eta}$ on the hyperplane $H_{\partial\Omega}({\bf y})$, which is not necessarily the projection $\Pi_{\partial\Omega}({\bf y})$. There is at least one orthogonal coordinate system centred on $\Pi_{\partial\Omega}({\bf y})$ in which ${\bf y}$ is expressed as ${\bf y}'=(-||\Pi_{\partial\Omega}({\bf y})-{\bf y}||,0,\ldots,0)$. Let ${\bm \eta}'$ be ${\bm \eta}$ in that frame. {\em If the particle is driven by a Brownian motion}, then the displacement along perpendicular directions to the first one is distributed according to \cite[section 3.3]{Buchmann_Petersen_SIAM}

 \begin{equation}
\label{F:II_Brownian_bridge}
\eta_j' \sim P_3= (y_j'-z_j')+\frac{\tau'}{h} \sqrt{\tau(1-\frac{\tau'}{h})}{\cal N} \qquad\textrm{for $j=2,...,D$},
\end{equation}

which is the density of a Brownian bridge pinned at $(y_j',0)$ and $(z_j',h)$ at some intermediate time $\tau'$. 
Algorithm \ref{A:BP} is the pseudocode for the Buchmann-Petersen method (assuming Brownian motion).  
In order to implement the rotations in an efficient way (which is especially critical in high dimensions), the reader is referred to \cite[section 3.4]{Buchmann_Petersen_SIAM}, where Givens rotations are used in order to construct and operate the orthogonal matrix $Q$. 

\begin{algorithm}[h!]
\caption{Method of Buchmann and Petersen for $\nabla^2u/2+f=0$ (BP)}
\begin{algorithmic}[1]
\STATE{{\bf Data:} $h>0$, ${\bf X}_0={\bf x}_0\in\Omega,{Y_0}=1,Z_0=0$, and $(d_0<0,{\bf N}_0)$}
\FOR{$k=0,1,2,\ldots$ until EXIT}
\STATE{Get ${\bf X}_{k+1},Y_{k+1}$ and $Z_{k+1}$ according to the Euler-Maruyama scheme (\ref{F:E_M})}
\STATE{Compute $(d_{k+1},{\bf N}_{k+1})$ according to the distance map of $\Omega$}
\IF{$d_{k+1}\geq 0$} 
	\STATE{Let $ \displaystyle{\gamma=\frac{d_k^2}h,\, \delta=\frac{|d_k|}{d_{k+1}}, \, w\sim {\cal IG}(\gamma,\delta),\textrm{  and } \tau_1=\frac{hw}{1+w}} $}
	\item[]
	\STATE{Let $\tau'=\tau_1$ and EXIT: {\bf goto} 15} 
\ELSE{} 
	\STATE{Let $w'\sim {\cal U}$ and $\displaystyle{\tau_2= \frac{-2|d_k|d_{k+1}}{\log{w'}}}$}
	\IF{$\tau_2<h$ (boundary test)}
		\STATE{Let $\tau'=\tau_2$ and EXIT: {\bf goto} 15}
	\ENDIF
\ENDIF
\ENDFOR
\item[] 
\STATE{Compute a $D\times D$ orthogonal matrix $Q$ such that $Q\big( \Pi_{\partial\Omega}({\bf X}_k)-{\bf X}_k \big)=(|d_k|,0,\ldots,0)^T$ using Givens rotations or Householder transformations \cite{Buchmann_Petersen_SIAM}}
\STATE{Let $(\eta'_1,\ldots,\eta'_D)={\bm \eta'}= Q({\bf X}_{k+1}-{\bf X}_k)$}
\FOR{j=2,\ldots,D}
	\STATE{Let $w''\sim {\cal N}$ and $\displaystyle{\eta''_j= \eta'_j + w''\sqrt{\tau'(1-\tau'/h)}}$}
\ENDFOR
\STATE{Let ${\bf X}_k'={\bf X}_k+Q^T{\bm \eta''}$ (such that ${\bf X}_k'\in H_{\partial\Omega}({\bf X}_k)$)}
\STATE{ Take $\nu=kh,{\bf X}_{\tau}=\Pi_{\partial\Omega}({\bf X}_k'),\,Y_{\tau}=Y_k,\,Z_{\tau}=Z_k$, and $\phi_h= g({\bf X}_{\nu})Y_{\nu}+Z_{\nu}$.}
\end{algorithmic}
\label{A:BP}
\end{algorithm}

The reasoning of Buchmann and Petersen can be extended to diagonal diffusions ({\em i.e.} in BVPs without crossed derivatives). Then, $\sigma$ is also diagonal and it is possible to {\em first} derive $\tau'$ from the component parallel to ${\bf N}_k$ from the analogous formula to (\ref{F:Tau1}) or (\ref{F:Tau2}), and {\em secondly} insert $\tau'$ in (\ref{F:II_Brownian_bridge}) in order to sample the remaining $D-1$ components on the tangent hyperplane. On the other hand, if the diffusion along the normal direction is not independent from the rest of orthogonal components (which is the general case, even with the frozen component approximation) the reasoning above fails. In theory, formulas (\ref{F:Tau1})-(\ref{F:II_Brownian_bridge}) might be extended to general diffusions based on high-dimensional Brownian bridges (with a full covariance matrix), which is a much more difficult problem. In any case, even in 1D, the method of Buchmann and Petersen with the frozen coefficient approximation may not improve on the Brownian bridge beyond the Brownian motion \cite{Buchmann_JCP}. 

\subsection{Milstein's walk on ellipsoids (WoE)}\label{SS:WoE}

Milstein's methods take bounded steps in order to avoid overshooting outside $\Omega$. Out of a variety of them (see \cite[chapter 6]{Milstein_Tretyakov_Book}), we have picked two. While both of them have been proved to yield a linear bias (under mild technical assumptions), this is dependent on the availability of geometrical and spectral information that may not be immediate. Perhaps for this reason, they are underreported--as far as we know, only \cite{Milstein_Tretyakov_Book} presents two numerical experiments. Here, we propose specific implementations of both, warning the reader that they rely on approximations which may spoil their theoretically linear weak convergence. In the WoE \cite[section 6.3.2]{Milstein_Tretyakov_Book}, the drift is first removed using Girsanov's theorem by taking ${\bm \mu}= \sigma^{-1}{\bm b}$ and ${\bm F}=0$ in (\ref{F:VR_Milstein_sys})--which may also affect the variance. Then, the system of SDEs to be solved reads:

\begin{equation}
\label{F:Girsanov_system}
\left\{
\begin{array}{ll}
d{\bf X}= \sigma({\bf X})d{\bf W}_t, & {\bf X}(0)={\bf x}_0, \\
dY= c({\bf X})Ydt + Y{\bf b}^T({\bf X})\sigma^{-T}({\bf X})d{\bf W}_t, & Y(0)=1,\\
dZ= f({\bf X})Ydt, & Z(0)= 0,
\end{array}
\right. 
\end{equation}

where $\sigma^{-T}=(\sigma^{-1})^T$. Let us define
\begin{equation}
\label{F:Elipsoide}
\theta({\bf x},{\bf y},s)=\{{\bf y}\,|\,
({\bf y}-{\bf x})^T\sigma^{-T}({\bf x})\sigma^{-1}({\bf x})({\bf y}-{\bf x})=s^2\}
=\{{\bf y}\,|\,
({\bf y}-{\bf x})^T A^{-1}({\bf x})({\bf y}-{\bf x})=2s^2\}.
\end{equation}

Since $A^{-1}$ is a positive-definite matrix, $\theta({\bf x},{\bf y},s)$ is a $D-$dimensional ellipsoid centred at ${\bf x}$ and with principal axes along the eigenvectors of $A^{-1}$. If $\lambda_{max}(A)\geq\lambda_2\geq\ldots \geq\lambda_D>0$ are the eigenvaues of $A$, the lengths of the semiaxes of $\theta({\bf x},{\bf y},s)$ are $s\sqrt{2\lambda_{max}}\geq\ldots \geq s\sqrt{2\lambda_D}>0$. Let us define the function ${\hat r}:\Omega\mapsto {\mathbb R}$ which determines at any point in the domain the tangent ellipsoid to $\partial\Omega$:

\begin{equation}
{\hat r}({\bf x})= \arg\min_s \theta({\bf x},{\bf y},s)\cap\partial\Omega\neq \emptyset.
\end{equation} 

Milstein's WoE has the sequence $\{{\bf X}_k\}_{k=0,1,...}$ hop over small ellipsoids of position-dependent size until coming close enough to the boundary, at which step $\partial\Omega$ is deemed hit--see Algorithm \ref{A:WoE}.

The discretization parameter in Algorithm \ref{A:WoE} is $r=\sqrt{Dh}$ rather than $h$. At any point ${\bf X}_k$ of the trajectory, the next hop is over the ellipsoid $\theta({\bf X}_k,{\bf y},r)$, except when ${\bf X}_k$ is so close to the boundary that $\theta({\bf X}_k,{\bf y},r)$ intersects the boundary, and the tangent ellipsoid $\theta({\bf X}_k,{\bf y},{\hat r}({\bf X}_k))$ has to be taken instead, in order to make sure that the trajectory will not overshoot the boundary. This happens when $|d_k|<r\sqrt{2\lambda_{max}(A)}$. (Note that $\lambda_d$ in \cite[chapter 6]{Milstein_Tretyakov_Book} is equivalent to $\sqrt{2\lambda_{max}(A)}$ here.) Eventually, the trajectory is stopped at the projection on the boundary when the distance is less than $r^2$ --see Figure \ref{I:Figura1} (right).  
 
Let ${\bm \omegaup}\sim{\cal B}(D)$ (look up Table \ref{T:Distributions}), and ${\bf y}=\sigma_k{\bm\omegaup}$. Since ${\bm\omegaup}^T{\bm\omegaup}=1={\bf y}^T\sigma_k^{-T}\sigma_k^{-1}{\bf y}$, ${\bf y}$ lies on the surface of the ellipsoid ${\bf y}^TA^{-1}{\bf y}=2$. Then, the possible new positions ${\bf X}_{k+1}={\bf X}_k+r_{k+1}\sigma({\bf X}_k){\bm \omegaup}_{k+1}$ in Algorithm \ref{A:WoE} are distributed over the surface of the ellipsoid $\theta({\bf X}_k,{\bf y},r_{k+1})$. Algorithm \ref{A:WoE} has an asymptotic weak convergence rate ${\cal O}(r^2)={\cal O}(h)$ for a $C^4$ solution of the BVP \cite[theorem 3.8]{Milstein_Tretyakov_Book}.

\begin{algorithm}[h!]
\caption{Milstein's walk on ellipsoids (WoE)}
\begin{algorithmic}[1]
\STATE{{\bf Data:} $h>0$, $r=\sqrt{Dh}$, ${\bf X}_0={\bf x}\in\Omega,{Y_0}=1,Z_0=0$, $(d_0<0,{\bf N}_0)$}
\FOR{$k=0,1,2,\ldots$ until EXIT}
\IF{$d_k\geq-r^2$}
\STATE{EXIT: {\bf goto} 15}
\ELSE \IF{$ -r\sqrt{2\lambda_{max}(A)} \leq d_k <-r^2$}
\STATE{Compute ${\hat r}({\bf X}_k)$ and let $r_{k+1}={\hat r}({\bf X}_k)$}
\ELSE{}
\STATE{$r_{k+1}=r$}
\ENDIF
\STATE{Let ${\bm \omegaup_{k+1}}\sim{\cal S}(D)$, and take one step inside $\Omega$ according to:}
\begin{equation}
\label{F:One_step_WoE}
\left\{
\begin{array}{l}
{\bf X}_{k+1}= {\bf X}_k + \sigma_k{\bm \omegaup_{k+1}}r_{k+1}, \\
Y_{k+1}= Y_k + Y_kc_k r_{k+1}^2/D +
Y_k{\bf b}_k^T\sigma^{-T}_k{\bm \omegaup_{k+1}}r_{k+1}, \\
Z_{k+1}= Z_k + Y_kf_k r_{k+1}^2/D. 
\end{array}
\right.
\end{equation}
\STATE{Compute $\{d_{k+1},{\bf N}_{k+1}\}$ according to the distance map of $\Omega$ and let $k=k+1$}
\ENDIF
\ENDFOR
\STATE{Let $\nu=(k+1)h, {\bf X}_{\nu}=\Pi_{\partial\Omega}({\bf X}_{k+1}), Y_{\nu}=Y_{k+1}, Z_{\nu}=Z_{k+1}$ and $\phi_h=g({\bf X}_{\nu})Y_{\nu}+Z_{\nu}$ } 
\end{algorithmic}
\label{A:WoE}
\end{algorithm}

In order to implement Algorithm \ref{A:WoE}, the calculation of ${\hat r}$ is necessary, and this is not clarified either in \cite{Milstein_Tretyakov_Book} or in \cite{Milstein97}. Since this geometric construction is difficult in the event of interaction with an irregularly shaped or nonsmooth boundary (and would likely be very time-consuming in that case), we propose here to take the half-space approximation (see the right side of Figure \ref{I:Figura1}, where for simplicity $\bf x=0$.) Let 
\begin{equation}
m({\bf y})= {\bf y}^T{\sigma}^{-T}\sigma^{-1}{\bf y}-r^2= {\bf y}^TA^{-1}{\bf y}-2r^2.
\end{equation}
The tangency point ${\bf y}_0$ belongs to the hyperplane which is at a distance ${\bf N}^T{\bf y}_0$ from the origin, and since $m({\bf y}_0)=0$
\begin{equation}
\label{F:elipf}
\nabla m({\bf y}_0)=2A^{-1}{\bf y}_0=2q {\bf N}
\end{equation}
for some $q\neq 0$. Since $\det\sigma^{-1}\neq 0$,
\begin{equation}
{\bf y}_0=\frac{k}{2}A{\bf N}\Rightarrow
\left\{
\begin{array}{l}
q^2{\bf N}^TA^{T}A^{-1}A{\bf N}= 2{\hat r}^2,
\textrm{ and }\\
{\bf N}^T{\bf y}_0=q{\bf N}^TA{\bf N},
\end{array}
\right.
\end{equation}
and recalling that ${\bf N}^T{\bf y}_0$, the distance to the tangent hyperplane, is the distance to the boundary in the half-space approximation,
\begin{equation}
{\hat r}({\bf X}_k)= \frac{|d({\bf X}_k)|}{\sqrt{2{\bf N}_k^TA^T{\bf N}_k}}=
\frac{|d({\bf X}_k)|}{||\sigma_k{\bf N}_k||}.
\end{equation}

Accurately determining $\lambda_{max}(A)$ may be time-consuming if $A=A({\bf x})$, specially in high dimensions. In those cases, a better approach may just be to bound it above using the Gershgorin's circles estimate (\ref{F:Gershgorin}). Even if $\lambda_{max}(A)$ is sharp, introducing the half-space approximation for the computation of $\theta({\bf X}_k,{\bf y},r_k)$ and thus ${\hat r}({\bf X}_k)$ in Algorithm \ref{A:WoE} may spoil the strict linear weak convergence rate of the WoE. The reason is that close to a corner or to a curved boundary, the tangent ellipsoid to the hyperplane tangent to boundary may stick slightly out of $\partial\Omega$, which might lead to the trajectory hopping outside $\Omega$. In that event, $d_{k+1}>0$ and line 3 in Algorithm \ref{A:WoE} terminates the trajectory.

As a last but important point, in (\ref{F:One_step_WoE}), the inverse matrix $\sigma({\bf X}_k)^{-T}$ must be calculated in order to compute $Y_{k+1}$. The practical approach is of course to solve the linear system ${\bf b}_k=\sigma_k{\bm \mu}_k$ with Gauss elimination, exploiting the fact that $\sigma$ can always be chosen triangular (Section \ref{SS:Gershgorin}). Still, this involves a cost of up to ${\cal O}(D^2)$ per step and trajectory, which may slow down the simulation if a high-dimensional BVP with many crossed derivatives is being solved.

\subsection{Milstein's random walk (RW)}\label{SS:RW}

Algorithm \ref{A:RW} summarizes algorithms 2.1 and 3.1 in \cite[chapter 6]{Milstein_Tretyakov_Book}. This method takes bounded steps whose components along the coordinate axes are quantized, for the components of the vector ${\bm \xi}$ in (\ref{F:Quantized_step}) are $\pm 1$ with equal probability $1/2$. In order to prevent the trajectory from taking a step outside $\Omega$, a {\em boundary zone} is defined as the set of points from which at least one of the $2^D$ possible steps defined by (\ref{F:Quantized_step}) would take the trajectory outside $\Omega$. Since this definition is non-constructive, one assumes that there is a positive scalar $\lambda$ such that, if ${\bf x}\in\Omega$ and $|d({\bf x})|\leq\lambda\sqrt{h}$, then ${\bf x}$ belongs to the boundary zone. Thus for practical purposes, the boundary zone is replaced by the interior shell $S_{\lambda}$ of thickness $\lambda\sqrt{h}$ which contains it. 

\begin{algorithm}[H]
\caption{Milstein's random walk (RW)}
\begin{algorithmic}[1]
\STATE{{\bf Data:} $h>0$, ${\lambda}>0$, ${\bf X}_0={\bf x}_0\in\Omega,{Y_0}=1,Z_0=0$, $(d_0<0,{\bf N}_0)$}
\FOR{$k=0,1,2,\ldots$ until EXIT}
\STATE{[SAFEGUARD: {\bf if} $d_k\geq 0$ {\bf then} EXIT: {\bf goto} 18. {\bf endif}]}
\IF{$d_k\geq -\lambda\sqrt{h}$} 
	\STATE{Let $w\sim{\cal U}$ and $p=\displaystyle{\frac{\lambda\sqrt{h}}{|d_k|+\lambda\sqrt{h}}}$}
	\IF{$w>p$}
		\STATE{Let ${\bf X}'_k= {\bf X}_k-\lambda\sqrt{h}{\bf N}_k$ (note that ${\bf X}'_k\in\Omega$ but ${\bf X}'_k\notin S_{\lambda}$)}
	\ELSE
	\STATE{The boundary is deemed hit at $\Pi_{\partial\Omega}({\bf X}_k)$. EXIT and {\bf goto} 18}
	\ENDIF
\ELSE
\STATE{Let ${\bf X}_k'={\bf X}_k$}
\ENDIF
\STATE{Let ${\bm\xi}=(\xi_1,\ldots,\xi_D)$ with $\xi_j\sim{\cal B}\,\, (j=1,\ldots,D)$}
\STATE{Take a step inside $\Omega$ according to:
\begin{eqnarray}\label{F:Quantized_step}
\left\{\begin{array}{l}
{\bf X}_{k+1}= {\bf X}'_k + h{\bf b}({\bf X}'_k) + \sqrt{h}\sigma({\bf X}'_k){\bm\xi},\\
Y_{k+1}= Y_k + hc({\bf X}'_k)Y_k,\\
Z_{k+1}= Z_k + hf({\bf X}'_k)Y_k.\\
\end{array}\right.
\end{eqnarray}}
\STATE{Compute $(d_{k+1},{\bf N}_{k+1})$ according to the distance map of $\Omega$}
\ENDFOR
\STATE{Let $\nu=kh, {\bf X}_{\nu}=\Pi_{\partial\Omega}({\bf X}_k), Y_{\nu}=Y_k, Z_{\nu}=Z_k$ and $\phi_h=g({\bf X}_{\nu})Y_{\nu}+Z_{\nu}$. } 
\end{algorithmic}
\label{A:RW}
\end{algorithm}

Whenever the trajectory enters $S_{\lambda}$, a special step is taken than either takes it on the closest boundary point, or throws it back inside $\Omega$ and out of $S_{\lambda}$. This binary event has an associated probability given by $p$ in Algorithm \ref{A:RW}, line 5. Eventually, the trajectory will hit the boundary after entering $S_{\lambda}$, and finish. This method has an asymptotic linear weak order of convergence \cite[section 6.3.1]{Milstein_Tretyakov_Book} under $C^2$  $\partial\Omega$ and $C^3$ SDE coefficients \cite[section 6.2.2]{Milstein_Tretyakov_Book}.

This method's snag is to determine $\lambda$, for the optimal way of doing it turns out to be a problem-dependent issue, not addressed in \cite{Milstein_Tretyakov_Book}. In order to reduce the bias, it is desirable that $\lambda$ be as small as possible as long as $S_{\lambda}$ contains the boundary zone. Note that
\begin{equation}
\label{F:RW_bound}
||\sigma{\bm\xi}||\leq ||{\bm \xi}||\sup_{{\bm\xi}\neq 0}\frac{||\sigma{\bm\xi}||}{||{\bm\xi}||}=
||{\bm \xi}||\sqrt{\lambda_{max}(\sigma^T\sigma)}=\sqrt{D\lambda_{max}(A)},  
\end{equation}
because $||{\bm \xi}||=\sqrt{D}$ and $\lambda_{max}(\sigma^T\sigma)=\lambda_{max}(\sigma\sigma^T)$. (The latter comes from the spectra of $\sigma^T\sigma$ and $\sigma\sigma^T$ being the same: if ${\bf v}$ is an eigenvector of $\sigma^T\sigma$ with eigenvalue $s$, then ${\bf w}=\sigma{\bf v}$ is an eigenvector of $\sigma\sigma^T$ with the same eigenvalue, and viceversa.) Using (\ref{F:RW_bound}), the steps can be bounded as
 \begin{equation}
||{\bf X}_{k+1}-{\bf X}_k||\leq
h||{\bf b}_k|| + \sqrt{Dh\lambda_{max}(A)}.
\end{equation}
Then, a spatially-dependent $\lambda=\lambda({\bf x})$ given by
\begin{equation}
\sqrt{h}\lambda({\bf x})= h||{\bf b}({\bf x})|| + \sqrt{Dh\lambda_{max}(A)} 
\end{equation}   
is an obvious choice. But this implicitly assumes that the ball centred at ${\bf x}$ with radius $\sqrt{h}\lambda({\bf x})$ is contained in $\Omega$--which is similar to a local half-space approximation. Moreover, an upper bound to $\lambda_{\max}(A)$ has to be estimated. To sum up, this approach amounts to replacing line 4 in Algorithm \ref{A:RW} by
\begin{equation}
\label{F:RW0}
d_k\geq -h||{\bf b}({\bf x})|| - \sqrt{Dh\lambda_{max}(A)}.
\end{equation}    
On the other hand, if a constant $\lambda$ is used, estimating it so that it is valid throughout $\Omega$ may be challenging, for it is an optimization problem involving the values that $\lambda_{max}\big(A({\bf x})\big)$ and the drift may take on inside a domain $\Omega$ with potentially irregular shape. Finally, as suggested in \cite{Milstein_Tretyakov_Book}, a last option is setting $\lambda$ to a number 'large enough', meaning that it will be valid on most steps, but allowing for the possibility that some trajectories get inadvertedly out of $\Omega$. Since that any of these options will most likely spoil the linear weak order of convergence, the safeguard in line $3$ of Algorithm \ref{A:RW} is required.

\section{Numerical experiments}\label{S:Experiments}

\subsection{Setup}\label{SS:Setup}

\begin{figure}[h]
\centerline{\includegraphics[width=.6\columnwidth]{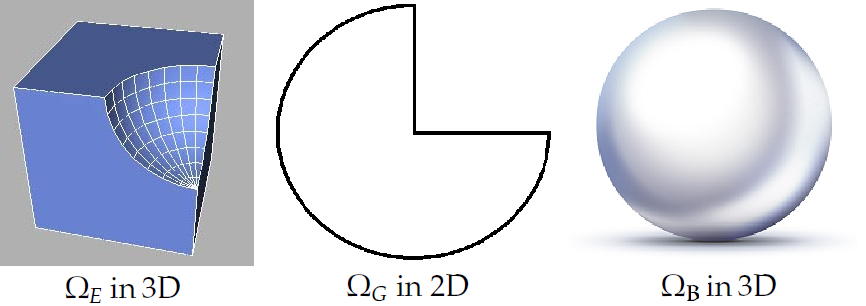}}
\caption{The three $D-$dimensional test domains. On the left is $\Omega_E$ in 3D. A ball with radius $\sqrt{D=3}/3$ and centred at $(1,1,1)$ has been removed, as well as another similar one at the origin (not visible).}
\label{I:Figura2}
\end{figure}

We consider three domains in ${\mathbb R}^D$, $D\geq 3$:
\begin{itemize}
\item $\Omega_B(R,{\bf C})$ is a hypersphere of radius $R$ and center ${\bf C}$.
\item $\Omega_G(R,{\bf C})$ is a hypersphere of radius $R$ at the origin from which the positive $D-$dimensional orthant has been carved out, and afterwards it has been shifted a vector ${\bf C}$, as in Figure \ref{I:Figura2} (in 2D). In 3D it looks like a piece of Gouda cheese.
\item $\Omega_E(L,{\bf C})$ looks in 3D like an piece of Emmental cheese (see Figure \ref{I:Figura2}), and is constructed from a hypercube with vertices at the origin and $(L,\ldots,L)$, where $L>0$. Then, two balls of radius one-third of the diagonal and centred on those vertices are removed, and the resulting body is shifted a vector ${\bf C}$.  
\end{itemize}
  
Note that while $\Omega_G$ has one single reentrant corner regardless of the dimension, $\Omega_E$ has ever more as $D$ increases. The code has been written in MATLAB in fully vectorized style, and run on single processors. It is made up of a main program, which manages the trajectories until exit and computes the Monte Carlo estimates; and a collection of files with the functions for the BVP to be solved, the domain, and the integrator. The main program input is a driver file with all the options for a given simulation. The whole suite will be made available from the first author's webpage.

The exact solution, $u_{ex}$, is always known and non-zero. The relative error is defined as
\begin{equation}
r_{h,N}= \frac{|u_{h,N}-u_{ex}|}{|u_{ex}|}.
\end{equation}

For every BVP and integrator we are interested in: 
\begin{enumerate}
\item The asymptotic weak order of convergence of $r_{h,N}$ as $h\shortrightarrow 0^+$ ({\em i.e.} $\delta$). In order to do this for each of the problems, we run $n$ independent simulations and construct an experimental plot of $r_{h,N}$ vs. the timestep $h=0.2,0.1,\ldots,0.2/2^{n-1}$. The simulations are run until the statistical error makes at most $20\%$ of the relative error with $95.5\%$ probability, {\em i.e.} $q\sqrt{Var[\phi_h({\bf x}_0)]/N}\leq |u_{h,N}-u_{ex}|/5$ with $q=2$. Then, $\delta$ is extracted from a loglog fit of those simulations for which $r_{h,N}<0.15$. Whenever bias cancellation obviously occurs (explained later in Section \ref{SS:Bias_Decomposition}), then the fit includes only the data points after the cancellation and with a relative error smaller than $15\%$. We always report the number of data points included in the fit along with the fitted value of $\delta$.    

\item The time it takes to produce a numerical solution $u_{h,N}$ with relative error $r_{h,N}$ less than $a$. We set $a$ depending on the BVP. Therefore, we choose $h$ from the previous convergence curve such that the bias is $a/2$, and $N$ such that the statistical error is smaller than $a/2$ with $95.5\%$ probability.    
\end{enumerate}

{\bf Variance reduction}. In order to speed up the simulations and being able to investigate the bias at small $h$, we apply variance reduction based on the exact formula (\ref{F:Pathwise_VR}). Specifically, we let
\begin{equation}
\label{F:Control_variates}
{\bm \mu}= 0,\, {\bf F}= -\sigma^T\nabla u_{ex},
\end{equation}

and modify the equation for $Z_{k+1}$ in each method accordingly. In the WoE, ${\bm\mu}\neq 0$, so that the formula that must be used is
\begin{equation}
\label{F:Combining_method}
{\bf F}= -\sigma^T\nabla u_{ex} - u_{ex}\sigma^{-1}{\bf b}.
\end{equation}

Note that in either case, the trajectories are the same as without variance reduction. While the SDE systems (\ref{F:Milstein_sys}) and (\ref{F:VR_Milstein_sys}) ({\em i.e.} with and without VR) are not identical, the numerical difference is, to all practical purposes, negligible, other than reducing (by up to $99.99\%$, or even more) the statistical error. We stress that variance reduction is used to construct the bias convergence curves, but obviously not when comparing the computational cost of the integrators for a given error tolerance $a$.

\subsection{Components of the discretization error (bias)}\label{SS:Bias_Decomposition}

Let us compactly write the BVP (\ref{F:EllipticBVP}) as
\begin{equation}
{\cal L}u+f=0 \textrm{ if }{\bf x}\in\Omega, \qquad u({\bf x})=g \textrm{ if }{\bf x}\in\partial\Omega,
\end{equation} 
where ${\cal L}:=\sum_{i=1}^D\sum_{j=1}^D a_{ij}\frac{\partial^2 u}{\partial x_i\partial x_j} + \sum_{k=1}^D b_k\frac{\partial u}{\partial x_k} + cu$. The solution $u$ can be split as

\begin{equation}
\label{F:Split}
u=v+w, \textrm{ with }
\left\{
\begin{array}{ll}
{\cal L}v+f=0, & v({\bf x}\in\partial\Omega)=0, \\
{\cal L}w=0, & w({\bf x}\in\partial\Omega)=g.
\end{array}
\right.
\end{equation}

From (\ref{F:gXYplusZ}), it can be seen that $g({\bf X}_{\tau})Y_{\tau}$ vanishes for $v$, while $Z_{\tau}=0$ for $w$. Therefore, when solving for $u$, $v$ and $w$ with a numerical integrator, the {\em signed} biases are also related as
\begin{equation}
\label{F:Bias_Split}
\lim_{N\shortrightarrow\infty}(u_{ex}-u_{h,N})= 
\lim_{N\shortrightarrow\infty}(v_{ex}-v_{h,N})+
\lim_{N\shortrightarrow\infty}(w_{ex}-w_{h,N}).
\end{equation} 
In other words, the bias of $u_{h,N}$ is the sum of two components: one influenced by the sampling of the Dirichlet BC by the trajectories (first term in (\ref{F:Dynkin})), and a second one which is the numerical approximation of the line integral for $f$ along the trajectory (second term in (\ref{F:Dynkin})). Both components, which we may call the BC sampling and quadrature errors, are respectively the signed biases of $w_{h,N}$ and $v_{h,N}$. Note that this is not to say that $v_{h,N}$ is unaffected by the presence of the boundary, for the length of individual trajectories is still determined by the first exit. However, since ${\cal L}$ is the same for $u,v$ and $w$, any $N$ trajectories involved in $u_{h,N},v_{h,N}$ and $w_{h,N}$ will be the same provided that the stream of random numbers is replicated upon computation of each of them. 

\begin{figure}[h]
\centerline{\includegraphics[width=.8\columnwidth]{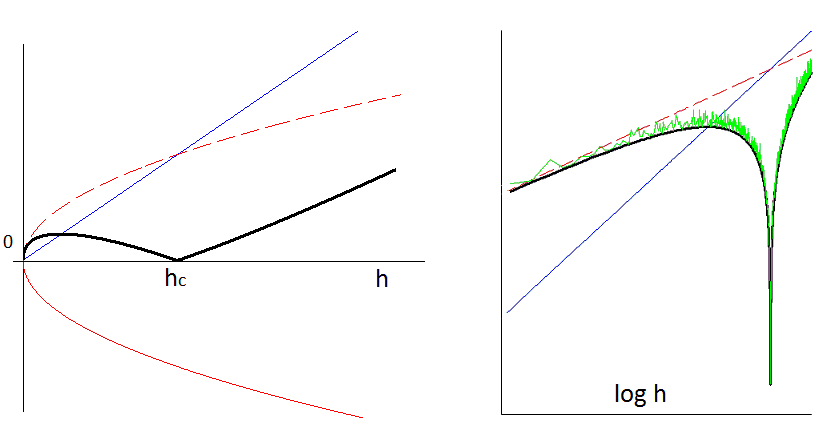}}
\caption{Bias decomposition and cancellation. Left: when the two components of the bias have opposed signs and exponents, a cancellation of the bias (thick line) may take place at finite $h=h_c$. Right: same plot in logarithmic scale. Note the 'plunge' and sharp 'rebound' of the bias. Gaussian noise proportional to one-fifth of the bias has been added (in green), to simulate the effect of the statistical error.}
\label{I:Figura3}
\end{figure}

This implicit decomposition of the bias often brings about a numerical artifact (so far unreported, to the best of our knowledge) that can lead to a misinterpretation of the behaviour of weak error curve with respect to $h$ when $h$ is not small enough; and which we explain next with an example. Let us denote the signed biases in (\ref{F:Bias_Split}) with ${\hat B}$, and assume that ${\hat B}(w_{h,N})=C_w\sqrt{h}$ and ${\hat B}(v_{h,N})=C_vh$. Then, $|{\hat B}(u_{h,N})|\shortrightarrow {\cal O}(\sqrt{h})$, but if the signs of $C_v$ and $C_w$ are opposed, then the bias of $u_{h,N}$ 'vanishes' at a finite value of $h=h_c>0$ where ${\hat B}(v_{h,N})$ and $-{\hat B}(w_{h,N})$ intersect. See Figure \ref{I:Figura3} for an illustration. In the loglog plot, the bias of $u_{h,N}$ (in absolute value) plunges, and then bounces back onto the slower rate of convergence (here, ${\cal O}(\sqrt{h})$). If data were only available for $h\gtrsim h_c$, the cancellation of the bias components might be mistaken for 'superlinear' convergence.           

The above artifact, which we will refer to as 'bias cancellation', accounts for the 'plunges' and 'rebounds' which can be observed in many of our experiments. Consider for instance the plunge of the error curve of the GM integrator at around $h\approx 10^{-4}$ in Figure \ref{I:Figura7}. 
In Table \ref{T:Bias_Decomposition} we indirectly demonstrate the artifact in that case. Let GM* be a variation of GM where, once the particle is deemed out at ${\bf X}_K$, the 'BC' $g$ is evaluated there rather than on $\Pi_{\partial\Omega}({\bf X}_K)$. (This is possible only because $g=u_{ex}|_{\partial\Omega}$, and $u_{ex}$ is available). The MC approximation yielded by GM*--let it be $u_{h,N}(GM*)$--only carries the quadrature and statistical errors, so that its signed bias is the same as that of $v_{h,N}(GM)$. As can be seen, the relative error of $u_{h,N}(GM*)$ is nonetheless larger than that of $u_{h,N}(GM)$ (as well as linearly decreasing). The reason why that the latter is smaller is the opposite contribution of the BC sampling error. The rightmost column in Table \ref{T:Bias_Decomposition} shows the difference (which is actually $u_{ex}-w_{h,N}(GM)$, plus some statistical noise), indicating that the BC sampling error starts dropping below $h\approx 10^{-4}$. Summing up, this worked out example illustrates that if various components of the bias converge at a different rate, cancellation, preceded by superlinearity, can take place.

\begin{table}[h!]
\begin{footnotesize}
\caption{Bias decomposition around the 'plunge' in the error convergence curve for the GM integrator on Figure \ref{I:Figura7}--see Section \ref{SS:Bias_Decomposition} for details.}
\label{T:Bias_Decomposition}
\end{footnotesize}
\[\begin{array}{|l|lll|}
\hline
h & u_{ex}-u_{h,N}(GM) & u_{ex}-u_{h,N}(GM*) & \textrm{difference} \\  
\hline
.000781   &  .00562\pm 4.5\times 10^{-5}            &  .00645\pm 4.3\times 10^{-5}  &  -.00083 \\
.000395  &  .00219\pm 3.2\times 10^{-5}            &  .00318\pm 3.0\times 10^{-5}  &  -.00010 \\
.000195  &  .00067\pm 2.3\times 10^{-5}            &  .00161\pm 2.1\times 10^{-5}  &  -.00095 \\
9.77\times 10^{-5}  &  4\times 10^{-6}\pm 1.6\times 10^{-5} &  .00080\pm 1.5\times 10^{-5}  &  -.00080 \\
4.88\times 10^{-5}  &  .00023\pm 1.2\times 10^{-5}            &  .00039\pm 1.1\times 10^{-5}  &  -.00016 \\ 
2.44\times 10^{-5} &  .00028\pm 8.2\times 10^{-6}            &  .00020\pm 7.5\times 10^{-6}  &   .00008 \\
\hline
\end{array}\]  
\end{table}

\subsection{Example I: General diffusion in 3D}\label{SS:ExampleI}
This is the test problem of Gobet and Menozzi in \cite{Gobet&Menozzi_2010}.It is a 3$D$ diffusion where all the coefficients are position dependent except for the potential $c$, which is zero. On it, Gobet and Menozzi observed linear convergence of their integrator. Like all examples, it will be solved on the three domains $\Omega_B,\Omega_E,\Omega_G$ (see Figure \ref{I:Figura2}). The diffusion matrix is given by
\begin{equation}
\label{F:sigma_Gobet_Menozzi}
\sigma= 
\left(
\begin{array}{ccc}
\sqrt{1+\mid z\mid} &  0  &  0 \\
\frac{1}{2}\sqrt{1+\mid x \mid} & \frac{\sqrt{3}}{2}\sqrt{1+\mid x \mid} & 0 \\
0  & \frac{1}{2}\sqrt{1+\mid y \mid}& \frac{\sqrt{3}}{2}\sqrt{1+\mid y \mid}   
\end{array}
\right), 
\end{equation}

while the rest of functions are:
\begin{equation}
{\bf b}=(y,z,x)^T,\,c=0,\,f= y^2z + z^2x + x^2y + \frac{1}{2}\sqrt{1+\mid z \mid}\sqrt{1+\mid x \mid}+\frac{\sqrt{3}}{2}x\sqrt{1+\mid x \mid}\sqrt{1+\mid y \mid}. 
\end{equation}

The solution is $u_{ex}=g=xyz$. 
Regarding Milstein's methods, $\sigma$ is already in triangular form as is convenient for computing ${\bm\mu}=\sigma^{-1}{\bf b}$. In any case, since $D=3$, $\sigma^{-1}$ and $\lambda_{max}(A)$ can be precalculated symbolically and stored, which will make the integration with WoE faster:
\begin{equation}
\sigma^{-1}=
\left(\begin{array}{ccc}
\frac{1}{\sigma_{11}} 							& 0 					& 0 			\\
\frac{-\sigma_{21}}{\sigma_{11}\sigma_{22}} 	& \frac{1}{\sigma_{22}} 	& 0 	\\
\frac{\sigma_{21}\sigma_{32}}{\sigma_{11}\sigma_{22}\sigma_{33}}     & \frac{-\sigma_{32}}{\sigma_{22}\sigma_{33}} & \frac{1}{\sigma_{33}}
\end{array}\right)
,\qquad
A=\frac{1}{2}
\left(\begin{array}{ccc}
\sigma_{11}^2 		& \sigma_{11}\sigma_{21} 		& 0 				\\
\sigma_{11}\sigma_{21}  & \sigma_{22}^2 + \sigma_{21}^2 	& \sigma_{22}\sigma_{32} 	\\
0			& \sigma_{22}\sigma_{32} 		& \sigma_{33}^2 + \sigma_{32}^2    
\end{array}\right).
\end{equation} 

We have solved this problem at ${\bf x}_0=(0.56,0.52,0.30)$ in $\Omega_B(1,[0,0,0])$; at ${\bf x}_0=(0.57,0.57,0.57)$ in $\Omega_G(1,[0.67,0.67,0.67])$; and at ${\bf x}_0=(0.57,0.57,0.57)$ in $\Omega_E(1,[0,0,0])$. Notice that the distances to the boundary are different in each case. 

\begin{figure}[h]
\centerline{\includegraphics[width=1\columnwidth]{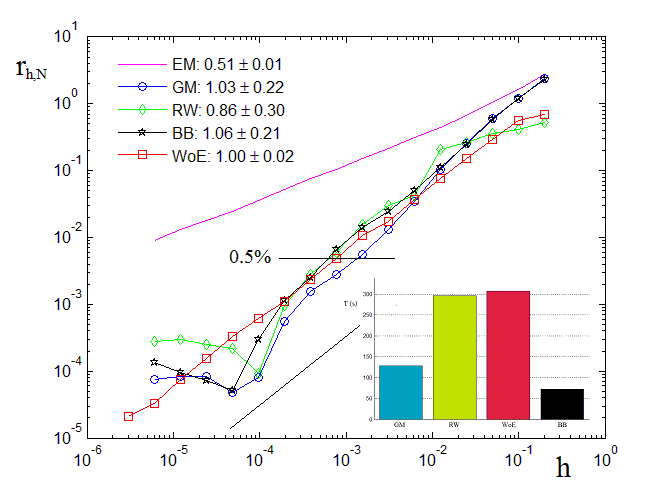}}
\caption{Relative error convergence, fitted $\delta$, and comparison of times with $a=0.01u_{ex}$ for Example I in $\Omega_B$. The segment of slope one is a guide for the eye.}
\label{I:Figura4}
\end{figure}

\begin{figure}[h]
\centerline{\includegraphics[width=1\columnwidth]{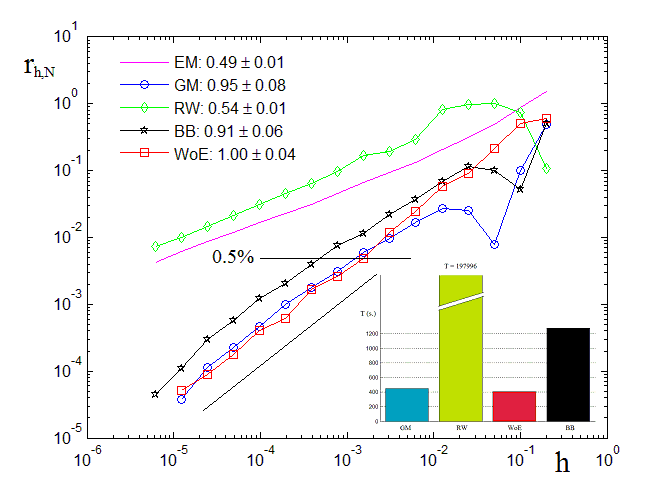}}
\caption{Relative error convergence, fitted $\delta$, and comparison of times with $a=0.01u_{ex}$ for Example I in $\Omega_G$.}
\label{I:Figura5}
\end{figure}

\begin{figure}[h]
\centerline{\includegraphics[width=1\columnwidth]{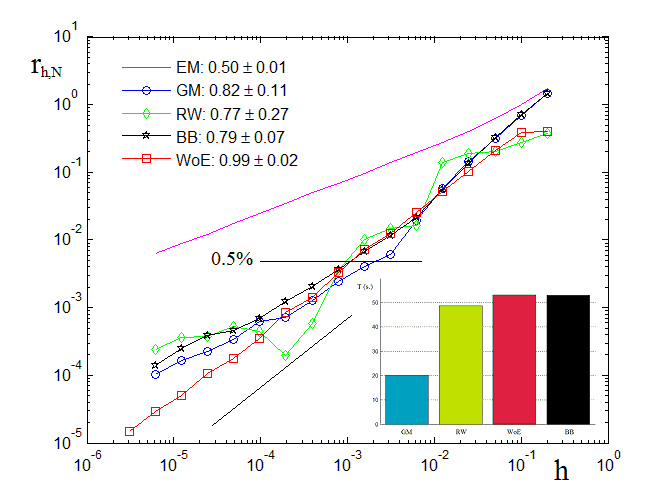}}
\caption{Relative error convergence, fitted $\delta$, and comparison of times with $a=0.01u_{ex}$ for Example I in $\Omega_E(1,[0,\ldots,0])$.}
\label{I:Figura6}
\end{figure}

Results are shown in Figures \ref{I:Figura4}-\ref{I:Figura6}. Error bars for $r_{h,N}$ are not shown for the sake of clarity, but recall that they are always $\pm r_{h,N}/5$. In this example, we have overlaid the convergence curves
with the estimates for $\delta$ (recall that not all the points are included into the $\delta$ fit), as well as with the computational times for an error tolerance $a=0.01u_{ex}$.

Here, RW uses the {\em ad-hoc} approximation $\lambda=2$. It works in $\Omega_B$ and $\Omega_E$, but definitely not in $\Omega_G$, due to the frequent overshoots of the inner corner. 

Several instances of bias cancellation are evident on those plots. In Figure \ref{I:Figura4}, all the integrators except EM and WoE 'bounce back' slightly below $h=10^{-4}$. In Figure \ref{I:Figura7}, the apparent 'superlinearity' of the WoE strongly hints at bias cancellation below $h=5\times 10^{-6}$. Unfortunately, it would take too long a simulation (well over $10^6$ seconds) to make sure.

\subsection{Example II: General diffusion in arbitrary dimension}\label{SS:ExampleII}

This is a synthetic problem with crossed second derivatives, oscillating coefficients and oscillating solution that can be easily extended to any arbitrary dimension. The exact solution is $u_{ex}=2+\cos{{\bf k}^T({\bf x}-{\bf y})}$. The evaluation point is ${\bf x}_0=(.7,.7,\ldots,.7)^T$, and the domains are $\Omega_G(1,[0.67,\ldots,0.67])$ and $\Omega_E(1,[0,\ldots,0])$. In order for the solution to be qualitatively similar regardless of $D$, we choose ${\bf k}=(1,2,3,1,2,\ldots,mod(D/3))^T$, where $mod(D/3)$ is the modulus, and we set ${\bf y}={\bf x}_0$, so that $u_{ex}({\bf x}_0)=3$. The diffusion matrix is chosen as
\begin{equation}
\label{F:Tril}
\sigma=\left(\begin{array}{cccc}
1 & 0  & \ldots & 0 \\
1 & 1  & \ldots & 0 \\
\vdots & \vdots & \ddots & \vdots \\
1 & 1  & \ldots & 1
\end{array}\right)
\Rightarrow
\sigma\sigma^T/2= A=\frac{1}{2} 
\left(\begin{array}{ccccc}
1 & 1 & 1 & \ldots & 1 \\
1 & 2 & 2 & \ldots & 2 \\
1 & 2 & 3 & \ldots & 3 \\
\vdots & \vdots & \vdots & \ddots & \vdots \\
1 & 2 & 3 & \ldots & D 
\end{array}\right).
\end{equation}
The rest of coefficients are
\begin{equation}
{\bf b}={\bf x}, \qquad c=-||{\bf x}||^2, 
\end{equation}
\begin{equation}
f= ({\bf k}^TA{\bf k}+||{\bf x}||^2)\cos{{\bf k}^T({\bf x}-{\bf y})}+2||{\bf x}||^2+{\bf k}^T{\bf x}\sin{{\bf k}^T({\bf x}-{\bf y})}.
\end{equation}

\begin{table}[h!]
\begin{footnotesize}
\caption{Fitted $\delta$ for Example II. The figure in parentheses is the number of data points included in the fit (see Section \ref{SS:Setup} for the criteria). A poor error bar on $\delta$ means that the
relative error had not entered the asymptotic power-law regime.} 
\label{T:Delta_ExampleII}
\[\begin{array}{|l|ccccc|}
\hline
\textrm{domain} & EM & GM & RW & BB & WoE \\  
\hline 
\textrm{$\Omega_E$ in 3D}& 0.55\pm .02 (11)& 1.01\pm .42 (12)& 0.83\pm .21 (13)& 0.88\pm .23 (12) & 0.98\pm .04 (7)\\
\textrm{$\Omega_E$ in 6D}& 0.57\pm .02 (9) & 1.28\pm .11 (10)& 0.69\pm .06 (10)& 1.19\pm .28 (10)& 0.98\pm .03 (10)  \\
\textrm{$\Omega_E$ in 12D}& 0.64\pm .08 (7) & 1.11\pm .08 (7) & 0.80\pm .06 (7)  & 1.25\pm .19 (7)& 1.00\pm .04 (9) \\
\textrm{$\Omega_E$ in 24D}& 0.22\pm .29 (10)& 0.94\pm .12 (6) & 1.15\pm .86 (6)  & 0.93\pm .09 (6) & 1.43\pm .10 (6)\\
\textrm{$\Omega_E$ in 48D}& 0.47\pm 1.47 (6)& 1.02\pm .17 (5) & 1.09\pm .25 (5)  & 0.77\pm .36 (5) & 0.88\pm .15 (7)\\
\hline 
\textrm{$\Omega_G$ in 3D}& 0.62\pm .05 (13)& 0.96\pm .05 (13)& 0.55\pm .02 (14))& 1.27\pm .02 (10) & 1.02 \pm .03 (10\\
\textrm{$\Omega_G$ in 6D}& 0.55\pm .29 (11)& 1.21\pm .13 (9) & 0.90\pm .33 (10)& 0.97\pm .24 (11) & 1.00 \pm .02 (11)\\
\textrm{$\Omega_G$ in 12D}& 0.44\pm .38 (9) & 1.12\pm .06 (8) & 0.71\pm .06 (8) ) & 1.18\pm .14 (7) & 1.11 \pm .05 (9) \\
\textrm{$\Omega_G$ in 24D}& \textrm{not enough data} & 1.42\pm .23 (5) & 0.67\pm .05 (4)  & 1.50\pm .48 (4) & 1.18 \pm .07 (7) \\
\textrm{$\Omega_G$ in 48D}& \textrm{not enough data} & 1.04\pm .10 (5) & 0.70\pm .49 (4) & 1.10\pm .15 (6)  & 1.08 \pm .16 (4)\\
\hline
\end{array}\]
\end{footnotesize}
\end{table}

\begin{figure}[h]
\centerline{\includegraphics[width=1\columnwidth]{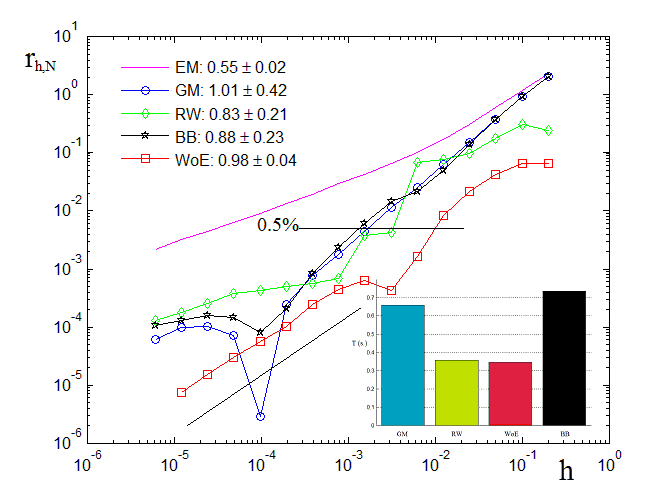}}
\caption{Relative error convergence, fitted $\delta$, and comparison of times with $a=0.01u_{ex}$ for Example II in $\Omega_E(1,[0,\ldots,0])$ (see Section \ref{SS:ExampleII}). The 'plunge' of the GM curve is worked out in Section \ref{SS:Bias_Decomposition}.}
\label{I:Figura7}
\end{figure}

\begin{figure}[h]
\centerline{\includegraphics[width=1.3\columnwidth]{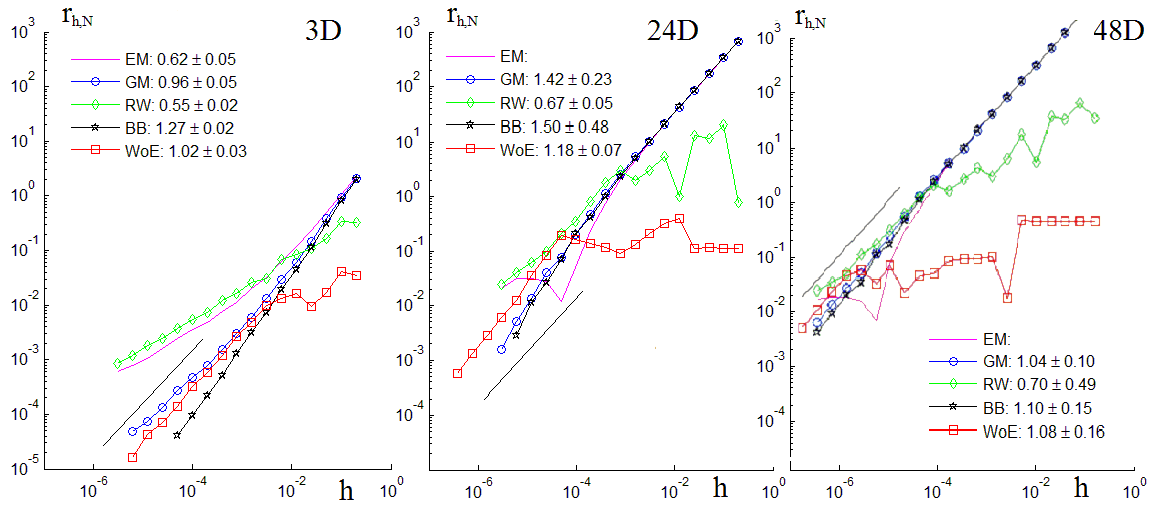}}
\caption{Convergence of the relative error and fitted $\delta$ for Example II in $D=\{3,24,48\}$. Notice that the scale is the same. There are not enough acceptable data to fit $\delta_{EM}$ in 24D and 48D.}
\label{I:Figura8}
\end{figure}

The convergence of relative errors for $\Omega_E$ in 3D is shown in Figure \ref{I:Figura7}, and those for $\Omega_G$ in $3D$, $24D$, and $48D$ are shown in Figure \ref{I:Figura8}. As $D$ grows, so does the bias (and the variance), and the fit of $\delta$ based on the acceptable available data becomes more unreliable. See Table \ref{T:Delta_ExampleII}.

Finally, we compare CPU times for $2\%$ relative accuracy in Table \ref{T:CPUt_ExampleII}. The fast growth of computational time (due to the larger bias and variance) along with increased dimension is also reported there.

\begin{table}[h!]
\begin{footnotesize}
\caption{Comparable computing times (in seconds) taken to solve Example II down to a relative statistical error equal or less than $1\%$. Same domains and ${\bf x}_0$ as in Table \ref{T:Delta_ExampleII}. WoE uses analytic formulas for $\sigma^{-1}$ and $\lambda_{max}(A)$, while WoE* solves the linear system with $\sigma$ and uses Gershgorin's circles.}
\label{T:CPUt_ExampleII}
\[\begin{array}{|l|lllll|l|}
\hline 
\textrm{Domain} & GM & RW & BB & WoE & WoE* \\  
\hline 
\textrm{$\Omega_E$ in $3D$} & 0.31 & 0.18  & 0.41 & 0.34 & 1.11 \\
\textrm{$\Omega_E$ in $6D$} & 6.8  & 12.9  & 6.1 & 9.58  & 14.03\\
\textrm{$\Omega_E$ in $12D$}& 153  & 296   & 121 & 216.50 & 327.85\\
\textrm{$\Omega_E$ in $24D$}& 3045 & 835   & 4392 & 3320.27 & 4149.98\\
\hline 
\textrm{s in: cost $\sim D^s$} & 4.3 & 1.7 & 5.0 & 4.9 & 3.7 \\
\hline
\textrm{$\Omega_G$ in $3D$} & 3.7  & 30.7  & 2.7 & 5.88 & 12.29 \\
\textrm{$\Omega_G$ in $6D$} & 53   & 31    & 31 & 93.89 & 128.62 \\
\textrm{$\Omega_G$ in $12D$}& 297  & 1708  & 341 & 596.77 & 853.59\\
\textrm{$\Omega_G$ in $24D$}& 1086 & \textrm{too long} & 1015 & 3199.99 & 4171.26\\
\hline 
\textrm{s in: cost $\sim D^s$} & 1.9 & - & 1.8 & 2.4 & 2.3 \\\hline 
\end{array}\]
\end{footnotesize}
\end{table}

Problem II suits itself well to assessing the importance of having analytical formulas for $\lambda_{max}(A)$ and $\sigma^{-1}$ when using the WoE in high dimensions. Here, WoE computes and stores both $\lambda_{\max}(A)$ and $\sigma^{-1}$ before start. Since $A$ in (\ref{F:Tril}) is constant, $\lambda_{\max}(A)$ can be precalculated numerically in any dimension, and the inverse of $\sigma$ in (\ref{F:Tril}) is trivial:

\begin{equation}
{\sigma}^{-1}=\left(\begin{array}{ccccc}
1  & 0  & 0 & \ldots & 0 \\
-1 & 1  & 0 & \ldots & 0 \\
0  & -1 & 1 & \ldots & 0 \\
\vdots  & \vdots & \ddots & \ddots & \vdots \\
0  & \ldots & \ldots & -1 & 1
\end{array}\right).
\end{equation}

On the other hand, in a BVP where $A=A({\bf x})$ instead of constant, there may be no analytical formulas available for $\lambda_{max}[A({\bf x})]$ and $\sigma^{-1}(\bf x)$. In that case, as pointed out in Section \ref{SS:WoE}, the best approach is to solve ${\bm b}(\bf x)=\sigma(\bf x){\bm \mu}(\bf x)$ (exploiting the row-echelon form of $\sigma$) and bound $\lambda_{max}[A({\bf x})]$ above using Gershgorin's theorem. In Table \ref{T:CPUt_ExampleII}, we have appended a column for the 'WoE*' integrator, which does this. Computational times for WoE* in Table \ref{T:CPUt_ExampleII} are definitely larger than for WoE, but only moderately so. This is encouraging, for it suggests that the WoE is viable also when analytical formulas for $\lambda_{max}[A({\bf x})]$ are not available.

\subsection{Examples III and IV: Poisson's equation in high dimension}\label{SS:ExampleIII}

Next, we tackle two problems solved by Buchmann and Petersen to test BP. First, we consider problem 4.2.1 in \cite{Buchmann_Petersen_SIAM} ({\bf Example III}):

\begin{equation}
\frac{1}{2}\nabla^2 u= -1, \qquad g=\sum_{i=1}^D x_i \Rightarrow u= (1-||{\bf x}||^2)/D + \sum_{i=1}^D x_i,
\end{equation}

in a $D-$dimensional ball $\Omega_B(1,[0,\ldots,0])$. Notice that $\sigma=I_D$ (so that $\lambda_{max}=1$), and that the exact solution is not an extension of the BC. Results for the weak order of convergence at ${\bf x}_0=(.9,0,\ldots,0)$ in $D=\{16,32,48,64\}$ are listed in Table \ref{T:Delta_ExampleIII}, and comparable computing times for $a=0.5\%$ in Table \ref{T:CPUt_ExampleIII}. RW and WoE use the exact value of $\lambda_{max}=\lambda(I)=1$. Here, computational times are only mildly sensitive to $D$ except for the WoE, even though the problem has no drift to be removed in the first place.

\begin{table}[h!]
\begin{footnotesize}
\caption{Fitted $\delta$ in Example III, and number of points included in the fit.}
\label{T:Delta_ExampleIII}
\[\begin{array}{|l|cccccc|}
\hline 
D & EM & GM & RW & WoE & BB & BP\\  
\hline 
16  & 0.52\pm.01 (17) & 1.06\pm.05 (14) & 0.53\pm .05 (14) & 1.00\pm .002 (13) & 0.99\pm.02 (14) & 1.11\pm .03 (13) \\
32  & 0.54\pm.02 (16) & 1.05\pm.03 (14) & 0.65\pm .04 (14) & 1.00\pm .002 (12) & 1.01\pm.01 (15) & 1.21\pm .06 (9)  \\
48  & 0.56\pm.02 (16) & 1.05\pm.04 (14) & 0.65\pm .04 (13) & 1.00\pm .001 (11) & 0.99\pm.02 (15) & 1.34\pm .06 (7)  \\
64  & 0.57\pm.02 (15) & 1.04\pm.03 (14) & 0.71\pm .07 (13) & 1.00\pm 9\times 10^{-4} (11) & 0.99\pm.01 (15) & 1.33\pm .08 (7)  \\
\hline 
\end{array}\]
\end{footnotesize}
\end{table}

\begin{table}[h!]
\begin{footnotesize}
\caption{Comparable computing times in seconds for solving Example III with $a=0.005u_{ex}$.}
\label{T:CPUt_ExampleIII}
\[\begin{array}{|l|llllll|}
\hline 
D & EM & GM & RW & WoE & BB & BP\\  
\hline 
16 & 500   &  11.6 &   103.1  &  94.3  &  25.6   &   19.3    \\ 
24 & 360   &  14.7 &    98.5  & 163.6  &  26.9   &   17.9    \\ 
32 & 300   &  15.5 &    83.8  & 291.3  &  42.7   &   5.9     \\ 
48 & 278   &  18.5 &    95.7  & 469.7  &  42.9   &   22.0    \\ 
\hline 
\end{array}\]
\end{footnotesize}
\end{table}

Finally, {\bf Example IV} is problem 4.2.2b in \cite{Buchmann_Petersen_SIAM}. The BC is a quartic polynomial, and again $\sigma=I_D$:

\begin{equation}
\frac{1}{2}\nabla^2 u= \frac{1}{6}\sum_{i=1}^{32} i(x_i)^2, \qquad g=u=\frac{1}{6}\sum_{i=1}^{32} i(x_i)^4. 
\end{equation}

\begin{figure}[h]
\centerline{\includegraphics[width=1\columnwidth]{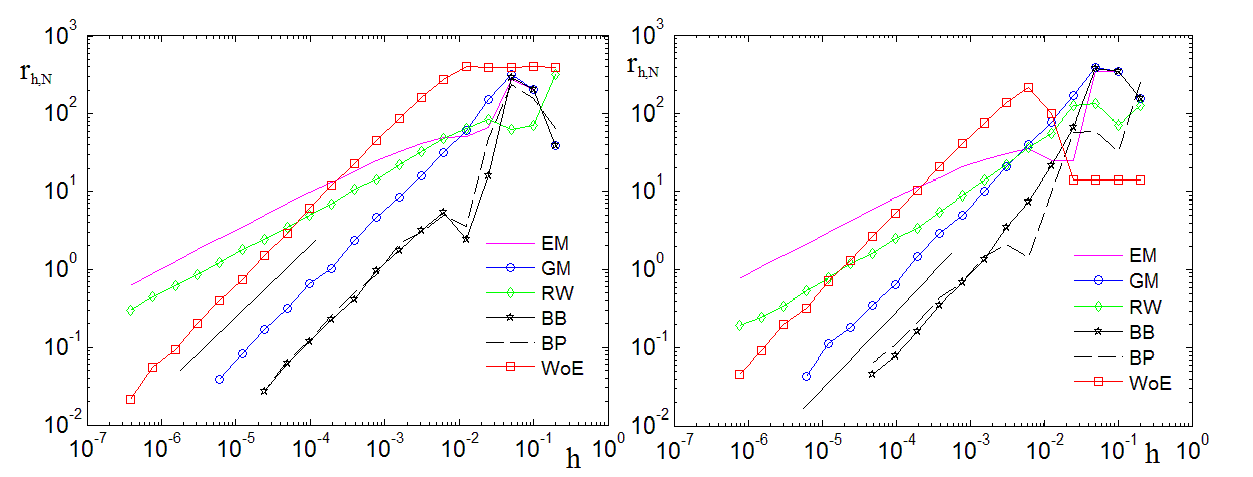}}
\caption{Convergence of relative errors for Problem IV in 32D $\Omega_B$ (left) and $\Omega_G$(right).}
\label{I:Figura9}
\end{figure}

\begin{figure}[hb!]
\centerline{\includegraphics[width=1\columnwidth]{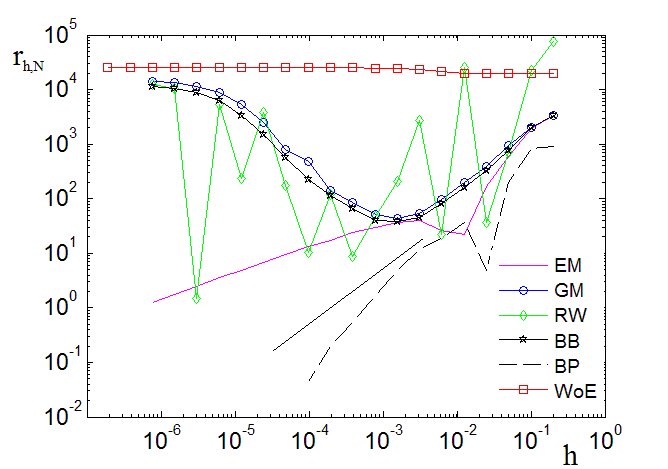}}
\caption{Convergence of relative errors for Problem IV in 32D $\Omega_E$. The point with the smallest error with BP took around 1.2 million seconds, the longest simulation reported in this paper.}
\label{I:Figura10}
\end{figure}

We take  ${\bf x}_0= (5/100,\ldots,5/100)\in {\mathbb R}^{32}$ in $\Omega_B(1,[0,\ldots,0])$, $\Omega_G(1,[1/10,\ldots,1/10])$, and $\Omega_E(1,[-1/2,\ldots,-1/2])$. (Approximate distances to the boundary are $0.72, 0.28, 0.45$, respectively.) The relative error behaviour is shown in Figures \ref{I:Figura9} and \ref{I:Figura10}. Here, we do not report $\delta$ because there are not enough data with $r_{h,N}\leq 0.15$ for a fit. For all the schemes but EM and RW, it looks like $\delta\shortrightarrow 1$.  While BP does not improve on BB in terms of accuracy in $\Omega_B$ and $\Omega_G$ (Figure \ref{I:Figura9}), we remark the striking results in $\Omega_E$ (Figure \ref{I:Figura10}). 
There, BP is the only integrator, along with EM, which converges at all. 

\section{Conclusions}\label{S:Conclusions}
   
In the light of the numerical evidence reported in this paper, it is difficult to pick an obvious winner from the five schemes compared. Instead, let us briefly comment on each of them:
\begin{itemize}
\item The boundary shift method by Gobet and Menozzi is reliable, fast and extremely easy to implement. 
Since it is basically Euler-Maruyama's method (on a slightly different domain), it inherits most of the available theory. For this reason, it was chosen for a Multilevel implementation \cite{Giles_y_yo}.

\item The Brownian bridge is also robust but somewhat less accurate (or slower) than the Gobet-Menozzi method, except if the diffusion is a Brownian motion. On the other hand, in the niche of Poisson equations, the method of Buchmann and Petersen is the scheme of choice, as clearly demonstrated by Example IV. Some numerical tests (not reported here) also suggest that BP yields important improvements in non-autonomous bounded SDEs. This is reasonable, due to the better accuracy in determining the exit time.  

\item Despite its simplicity, the Random Walk method has been a disappointment, for it seems to be very sensitive to parameters not always easy to approximate and to the presence of corners in the domain.  

\item Finally, Milstein's Walk on Ellipsoids--here combined with our half-space approximation to determine the largest tangent ellipsoid--seems to us the most stable in terms of yielding $\delta=1$ in nearly all the experiments, regardless of the dimension and the geometry. 
On the other hand, for an SDE involving a complicated diffusion and geometry (such as Examples I and II) where a really small tolerance $a$ is required, WoE would be our choice. Moreover, it is the basis for a linear scheme for general SDEs reflected on the boundary \cite[Section 6.6]{Milstein_Tretyakov_Book}.    
\end{itemize}

The point of this paper has been to show that weak orders of convergence better than ${\cal O}(\sqrt{h})$ for numerical approximations of bounded stopped diffusions are in general attainable, and quite often linear (or nearly so). This translates into large savings in computational times. The most significant challenge now is to raise the weak order of convergence in SDEs featuring non-smooth coefficients, reflecting BCs, or boundaries with sharp corners. Another important direction for improvement is the analysis of the the strong convergence properties of higher-order integrators, in order to insert them into the Multilevel framework.

\section{Acknowledgements}
This work was supported by Portuguese national funds through FCT under grants UID/CEC/50021/2013 and PTDC/EIA-CCO/098910/2008. FB also acknowledges FCT funding under grant SFRH/BPD/79986/2011.

 

\begin{thebibliography}{10}

\bibitem{Acebron2005}
J.A. Acebr\'on, M.P. Busico, P. Lanucara, and R. Spigler, {\em Domain decomposition solution of elliptic boundary-value problems via Monte Carlo and quasi-Monte Carlo methods}. SIAM J. Sci. Comp. {\bf 27}(2), 440-457 (2005). 

\bibitem{Baldi1995}
P. Baldi, {\em Exact asymptotics for the probability of exit from a domain and applications to simulation}. The Annals of Probability {\bf 23}(4), 1644-1670 (1995).

\bibitem{Bayer2010}
C. Bayer, A. Szepessy and R. Tempone, {\em Adaptive weak approximation of reflected and stopped diffusions}. Monte Carlo Methods and Applications {\bf 16}, 1-67 (2010).

\bibitem{Bernal2014}
F. Bernal, J.A. Acebr\'on and I. Anjam, {\em A stochastic algorithm based on Fast Marching for automatic capacitance extraction in non-Manhattan geometries}.
SIAM Journal on Imaging Sciences {\bf 7}(4), 2657-2674 (2014).

\bibitem{Buchmann_Petersen_SIAM}
F.M. Buchmann and W.P. Petersen, {\em An exit probability approach to solving high dimensional Dirichlet problems}. SIAM J. Sci. Comput. {\bf 28}, 1153-1166 (2006). 

\bibitem{Buchmann_JCP}
F.M. Buchmann, {\em Simulation of stopped diffusions}. Journal of Computational Physics {\bf 202}(2), 446-462 (2005).

\bibitem{Constantini1998}
C. Constantini, B. Pacchiarotti and F. Sartoretto, {\em Numerical approximation for functionals of reflecting diffusion processes}. SIAM J. Appl. Math. {\bf 58}, 73-102 (1998).

\bibitem{Deaconu2006}
M. Deaconu and A. Lejay, {\em A Random Walk on Rectangles algorithm}, Methodol. Comput. Appl. Probab. {\bf 8}, 135-151 (2006).

\bibitem{Freidlin_Book}
M. Freidlin, Functional Integration and Partial Differential Equations. 
Princeton University Press (1985).

\bibitem{Giles_seminal}
M.B. Giles, {\em Multi-level Monte Carlo path simulation}. Operations Research, {\bf 56}(3), 607-617 (2008).

\bibitem{Giles_y_yo}
M.B. Giles and F. Bernal, {\em Multilevel simulations of expected exit times and other functionals of stopped diffusions}. In preparation.

\bibitem{Giraudo1999}
M.T. Giraudo, L. Sacerdote, {\em An improved technique for the simulation of first passage times for diffusion processes}. Comm.
Statist. Simulation Comput. {\bf 28}(4), 1135-1163 (1999).

\bibitem{Glasserman2003}
P. Glasserman, Monte Carlo Methods in Financial Engineering. Springer (2003).

\bibitem{Gobet&Menozzi_2010}
E. Gobet and S. Menozzi, {\em Stopped diffusion processes: overshoots and boundary correction}. Stochastic Processes and their Applications, {\bf 120}, 130-162, (2010).

\bibitem{Gobet2001_killing}
E. Gobet, {\em Euler schemes for the weak approximation of killed diffusion}. Stochastic Process. Appl. {\bf 87}, 167-197 (2000).

\bibitem{Gobet2001_BB}
E. Gobet, {\em Euler schemes and half-space approximation for the simulation of diffusions in a domain}. ESAIM: Probability and Statistics {\bf 5}, 261-297 (2001).


\bibitem{Higham2013}
D.J. Higham, X. Mao, M. Roj, Q. Song, and G. Yin. {\em Mean exit times and the multilevel Monte Carlo method}. SIAM Journal on Uncertainty Quantification, {\bf 1}(1) 2-18 (2013).

\bibitem{Iyer2010}
G. Iyer, A. Novikov, L. Ryzhik, A. Zlatos, {\em Exit times for diffusions with incompressible drift}. SIAM J. Math. Anal. {\bf 42}, 2484-2498 (2010).

\bibitem{Jansons_Lythe_1st}
K.M. Jansons, G.D. Lythe, {\em Efficient numerical solution of stochastic differential equations using exponential timestepping} J. Statist. Phys. {\bf 100}(5-6) 1097-1109 (2000).

\bibitem{Kloeden&Platten}
P.E. Kloeden and E.Platten, Numerical Solution of Stochastic Differential Equations. Springer (1999).

\bibitem{Lemaire&Pages}
P. Lemaire and G. Pag\'es, {\em Multilevel Richardson-Romberg extrapolation}, to appear in Bernoulli, ArXiv 1401.1177v3 (2013).


\bibitem{Maire2008}
S. Mair\'e, E. Tanr\'e, {\em Some new simulation schemes for the evaluation of Feynman-Kac representations}. Monte Carlo Methods and Applications, {\bf 14}(1), 29-51 (2008).

\bibitem{Maire&Simon}
S. Mair\'e and M. Simon, {\em A partially reflecting random walk on spheres algorithm for electrical impedance tomography}, preprint in arXiv:1502.04318 (2015).

\bibitem{LTMCR}
S. Mancini, F. Bernal and J.A. Acebr\'on, {\em An efficient algorithm for accelerating Monte Carlo approximations of the solution to boundary value
problems}. Journal of Scientific Computing {\bf 66}(2) 577-597 (2016).


\bibitem{Mannella1999}
R. Mannella, {\em Absorbing boundaries and optimal stopping in a stochastic differential equation}. Physics Letters A {\bf 254}(5), 257-262 (1999). 

\bibitem{Mascagni2002}
M. Mascagni and N.A. Simonov, {\em Monte Carlo methods for calculating
some physical properties of large molecules}, SIAM J. Sci. Comput. {\bf 26}
, 339-357 (2004).

\bibitem{Michael1976}                                                                    
J.R. Michael, W.R. Schucany and R.W. Haas, {\em Generating random variates using transformations with multiple roots}. Am. Stat. {\bf 30}, 88-90 (1976).

\bibitem{Milstein97} 
G.N. Milstein, {\em Weak approximation of a diffusion process in a bounded domain}. Stoch. Stoch. Rep. {\bf 62} 147-200 (1997).

\bibitem{Milstein_Tretyakov_Book}
G.N. Milstein and M.V. Tretyakov, Stochastic Numerics for Mathematical Physics. Springer, Berlin (2004). 

\bibitem{Milstein&Tretyakov_simplest}
G.N. Milstein and M.V. Tretyakov, {\em The simplest random walks for the Dirichlet problem}. Teor. Veroyatnost. i Primenen. {\bf 47}(1), 39-58 (2002).

\bibitem{Miranda70} 
C. Miranda, Partial Differential Equations of Elliptic Type. Springer (1970).


\bibitem{Muller1956}
M.E. Muller, {\em Some continuous Monte Carlo methods for the Dirichlet
problem}. Ann. Math. Statist. {\bf 27}, 569-589 (1956).

\bibitem{Redner2001}
S. Redner, A Guide to First-Passage Processes. Cambridge University Press (2001).

\bibitem{Schwabedal_PRL}
J.T.C. Schwabedal and A. Pikovsky, {\em Phase description of stochastic oscillations}, Phys. Rev. Lett. {\bf 110}(20):204102 (2013).

\bibitem{Tamborrino2014}
M. Tamborrino, L. Sacerdote, and M. Jacobsen, {\em Weak
convergence of marked point processes generated by crossings of
multivariate jump processes. Application to neural network modeling},
Physica D, {\bf 288}, 45-52 (2014).

\bibitem{Talay&Tubaro90}
D. Talay and L. Tubaro, {\em Expansion of the global error for numerical
schemes solving stochastic differential equations}. Stoch. Anal. and App.,
{\bf 8}(4), 94-120 (1990).


\end{thebibliography}
\end{document}